\documentclass{amsart}

\numberwithin{equation}{section}

\usepackage{graphicx,amsthm}
\usepackage{verbatim,amsmath,amscd,amssymb,color}
\usepackage[mathscr]{eucal}
\input xy
\xyoption{all}
\begin{document}
\renewcommand{\labelenumi}{$($\roman{enumi}$)$}
\renewcommand{\labelenumii}{$(${\rm \alph{enumii}}$)$}
\font\germ=eufm10
\newcommand{\cI}{{\mathcal I}}
\newcommand{\cA}{{\mathcal A}}
\newcommand{\cB}{{\mathcal B}}
\newcommand{\cC}{{\mathcal C}}
\newcommand{\cD}{{\mathcal D}}
\newcommand{\cE}{{\mathcal E}}
\newcommand{\cF}{{\mathcal F}}
\newcommand{\cG}{{\mathcal G}}
\newcommand{\cH}{{\mathcal H}}
\newcommand{\cK}{{\mathcal K}}
\newcommand{\cL}{{\mathcal L}}
\newcommand{\cM}{{\mathcal M}}
\newcommand{\cN}{{\mathcal N}}
\newcommand{\cO}{{\mathcal O}}
\newcommand{\cR}{{\mathcal R}}
\newcommand{\cS}{{\mathcal S}}
\newcommand{\cV}{{\mathcal V}}
\newcommand{\cX}{{\mathcal X}}
\newcommand{\fra}{\mathfrak a}
\newcommand{\frb}{\mathfrak b}
\newcommand{\frc}{\mathfrak c}
\newcommand{\frd}{\mathfrak d}
\newcommand{\fre}{\mathfrak e}
\newcommand{\frf}{\mathfrak f}
\newcommand{\frg}{\mathfrak g}
\newcommand{\frh}{\mathfrak h}
\newcommand{\fri}{\mathfrak i}
\newcommand{\frj}{\mathfrak j}
\newcommand{\frk}{\mathfrak k}
\newcommand{\frI}{\mathfrak I}
\newcommand{\fm}{\mathfrak m}
\newcommand{\frn}{\mathfrak n}
\newcommand{\frp}{\mathfrak p}
\newcommand{\fq}{\mathfrak q}
\newcommand{\frr}{\mathfrak r}
\newcommand{\frs}{\mathfrak s}
\newcommand{\frt}{\mathfrak t}
\newcommand{\fru}{\mathfrak u}
\newcommand{\frA}{\mathfrak A}
\newcommand{\frB}{\mathfrak B}
\newcommand{\frF}{\mathfrak F}
\newcommand{\frG}{\mathfrak G}
\newcommand{\frH}{\mathfrak H}
\newcommand{\frJ}{\mathfrak J}
\newcommand{\frN}{\mathfrak N}
\newcommand{\frP}{\mathfrak P}
\newcommand{\frT}{\mathfrak T}
\newcommand{\frU}{\mathfrak U}
\newcommand{\frV}{\mathfrak V}
\newcommand{\frX}{\mathfrak X}
\newcommand{\frY}{\mathfrak Y}
\newcommand{\frZ}{\mathfrak Z}
\newcommand{\rA}{\mathrm{A}}
\newcommand{\rC}{\mathrm{C}}
\newcommand{\rd}{\mathrm{d}}
\newcommand{\rB}{\mathrm{B}}
\newcommand{\rD}{\mathrm{D}}
\newcommand{\rE}{\mathrm{E}}
\newcommand{\rH}{\mathrm{H}}
\newcommand{\rK}{\mathrm{K}}
\newcommand{\rL}{\mathrm{L}}
\newcommand{\rM}{\mathrm{M}}
\newcommand{\rN}{\mathrm{N}}
\newcommand{\rR}{\mathrm{R}}
\newcommand{\rT}{\mathrm{T}}
\newcommand{\rZ}{\mathrm{Z}}
\newcommand{\bbA}{\mathbb A}
\newcommand{\bbB}{\mathbb B}
\newcommand{\bbC}{\mathbb C}
\newcommand{\bbG}{\mathbb G}
\newcommand{\bbF}{\mathbb F}
\newcommand{\bbH}{\mathbb H}
\newcommand{\bbP}{\mathbb P}
\newcommand{\bbN}{\mathbb N}
\newcommand{\bbQ}{\mathbb Q}
\newcommand{\bbR}{\mathbb R}
\newcommand{\bbV}{\mathbb V}
\newcommand{\bbZ}{\mathbb Z}
\newcommand{\adj}{\operatorname{adj}}
\newcommand{\Ad}{\mathrm{Ad}}
\newcommand{\Ann}{\mathrm{Ann}}
\newcommand{\rcris}{\mathrm{cris}}
\newcommand{\ch}{\mathrm{ch}}
\newcommand{\coker}{\mathrm{coker}}
\newcommand{\diag}{\mathrm{diag}}
\newcommand{\Diff}{\mathrm{Diff}}
\newcommand{\Dist}{\mathrm{Dist}}
\newcommand{\rDR}{\mathrm{DR}}
\newcommand{\ev}{\mathrm{ev}}
\newcommand{\Ext}{\mathrm{Ext}}
\newcommand{\cExt}{\mathcal{E}xt}
\newcommand{\fin}{\mathrm{fin}}
\newcommand{\Frac}{\mathrm{Frac}}
\newcommand{\GL}{\mathrm{GL}}
\newcommand{\Hom}{\mathrm{Hom}}
\newcommand{\hd}{\mathrm{hd}}
\newcommand{\rht}{\mathrm{ht}}
\newcommand{\id}{\mathrm{id}}
\newcommand{\im}{\mathrm{im}}
\newcommand{\inc}{\mathrm{inc}}
\newcommand{\ind}{\mathrm{ind}}
\newcommand{\coind}{\mathrm{coind}}
\newcommand{\Lie}{\mathrm{Lie}}
\newcommand{\Max}{\mathrm{Max}}
\newcommand{\mult}{\mathrm{mult}}
\newcommand{\op}{\mathrm{op}}
\newcommand{\ord}{\mathrm{ord}}
\newcommand{\pt}{\mathrm{pt}}
\newcommand{\qt}{\mathrm{qt}}
\newcommand{\rad}{\mathrm{rad}}
\newcommand{\res}{\mathrm{res}}
\newcommand{\rgt}{\mathrm{rgt}}
\newcommand{\rk}{\mathrm{rk}}
\newcommand{\SL}{\mathrm{SL}}
\newcommand{\soc}{\mathrm{soc}}
\newcommand{\Spec}{\mathrm{Spec}}
\newcommand{\St}{\mathrm{St}}
\newcommand{\supp}{\mathrm{supp}}
\newcommand{\Tor}{\mathrm{Tor}}
\newcommand{\Tr}{\mathrm{Tr}}
\newcommand{\wt}{\mathrm{wt}}
\newcommand{\Ab}{\mathbf{Ab}}
\newcommand{\Alg}{\mathbf{Alg}}
\newcommand{\Grp}{\mathbf{Grp}}
\newcommand{\Mod}{\mathbf{Mod}}
\newcommand{\Sch}{\mathbf{Sch}}\newcommand{\bfmod}{{\bf mod}}
\newcommand{\Qc}{\mathbf{Qc}}
\newcommand{\Rng}{\mathbf{Rng}}
\newcommand{\Top}{\mathbf{Top}}
\newcommand{\Var}{\mathbf{Var}}
\newcommand{\gromega}{\langle\omega\rangle}
\newcommand{\lbr}{\begin{bmatrix}}
\newcommand{\rbr}{\end{bmatrix}}
\newcommand{\forb}{\bigcirc\kern-2.8ex \because}
\newcommand{\forbb}{\bigcirc\kern-3.0ex \because}
\newcommand{\forbbb}{\bigcirc\kern-3.1ex \because}
\newcommand{\cd}{commutative diagram }
\newcommand{\SpS}{spectral sequence}
\newcommand\C{\mathbb C}
\newcommand\hh{{\hat{H}}}
\newcommand\eh{{\hat{E}}}
\newcommand\F{\mathbb F}
\newcommand\fh{{\hat{F}}}
\newcommand\Z{{\mathbb Z}}
\newcommand\Zn{\Z_{\geq0}}
\newcommand\et[1]{\tilde{e}_{#1}}
\newcommand\ft[1]{\tilde{f}_{#1}}

\def\ge{\frg}
\def\AA{{\mathcal A}}
\def\al{\alpha}
\def\bq{B_q(\ge)}
\def\bqm{B_q^-(\ge)}
\def\bqz{B_q^0(\ge)}
\def\bqp{B_q^+(\ge)}
\def\beneme{\begin{enumerate}}
\def\beq{\begin{equation}}
\def\beqn{\begin{eqnarray}}
\def\beqnn{\begin{eqnarray*}}
\def\bigsl{{\hbox{\fontD \char'54}}}
\def\bbra#1,#2,#3{\left\{\begin{array}{c}\hspace{-5pt}
#1;#2\\ \hspace{-5pt}#3\end{array}\hspace{-5pt}\right\}}
\def\cd{\cdots}
\def\CC{\mathbb{C}}
\def\CBL{\cB_L(\TY(B,1,n+1))}
\def\CBM{\cB_M(\TY(B,1,n+1))}
\def\CVL{\cV_L(\TY(D,1,n+1))}
\def\CVM{\cV_M(\TY(D,1,n+1))}
\def\ddd{\hbox{\germ D}}
\def\del{\delta}
\def\Del{\Delta}
\def\Delr{\Delta^{(r)}}
\def\Dell{\Delta^{(l)}}
\def\Delb{\Delta^{(b)}}
\def\Deli{\Delta^{(i)}}
\def\Delre{\Delta^{\rm re}}
\def\ei{e_i}
\def\eit{\tilde{e}_i}
\def\eneme{\end{enumerate}}
\def\ep{\epsilon}
\def\eeq{\end{equation}}
\def\eeqn{\end{eqnarray}}
\def\eeqnn{\end{eqnarray*}}
\def\fit{\tilde{f}_i}
\def\FF{{\rm F}}
\def\ft{\tilde{f}_}
\def\gau#1,#2{\left[\begin{array}{c}\hspace{-5pt}#1\\
\hspace{-5pt}#2\end{array}\hspace{-5pt}\right]}
\def\gl{\hbox{\germ gl}}
\def\hom{{\hbox{Hom}}}
\def\ify{\infty}
\def\io{\iota}
\def\kp{k^{(+)}}
\def\km{k^{(-)}}
\def\llra{\relbar\joinrel\relbar\joinrel\relbar\joinrel\rightarrow}
\def\lan{\langle}
\def\lar{\longrightarrow}
\def\max{{\rm max}}
\def\lm{\lambda}
\def\Lm{\Lambda}
\def\mapright#1{\smash{\mathop{\longrightarrow}\limits^{#1}}}
\def\Mapright#1{\smash{\mathop{\Longrightarrow}\limits^{#1}}}
\def\mm{{\bf{\rm m}}}
\def\nd{\noindent}
\def\nn{\nonumber}
\def\nnn{\hbox{\germ n}}
\def\catob{{\mathcal O}(B)}
\def\oint{{\mathcal O}_{\rm int}(\ge)}
\def\ot{\otimes}
\def\op{\oplus}
\def\opi{\ovl\pi_{\lm}}
\def\osigma{\ovl\sigma}
\def\ovl{\overline}
\def\plm{\Psi^{(\lm)}_{\io}}
\def\qq{\qquad}
\def\q{\quad}
\def\qed{\hfill\framebox[2mm]{}}
\def\QQ{\mathbb Q}
\def\qi{q_i}
\def\qii{q_i^{-1}}
\def\ra{\rightarrow}
\def\ran{\rangle}
\def\rlm{r_{\lm}}
\def\ssl{\hbox{\germ sl}}
\def\slh{\widehat{\ssl_2}}
\def\ti{t_i}
\def\tii{t_i^{-1}}
\def\til{\tilde}
\def\tm{\times}
\def\tt{\frt}
\def\TY(#1,#2,#3){#1^{(#2)}_{#3}}
\def\ua{U_{\AA}}
\def\ue{U_{\vep}}
\def\uq{U_q(\ge)}
\def\uqp{U'_q(\ge)}
\def\ufin{U^{\rm fin}_{\vep}}
\def\ufinp{(U^{\rm fin}_{\vep})^+}
\def\ufinm{(U^{\rm fin}_{\vep})^-}
\def\ufinz{(U^{\rm fin}_{\vep})^0}
\def\uqm{U^-_q(\ge)}
\def\uqmq{{U^-_q(\ge)}_{\bf Q}}
\def\uqpm{U^{\pm}_q(\ge)}
\def\uqq{U_{\bf Q}^-(\ge)}
\def\uqz{U^-_{\bf Z}(\ge)}
\def\ures{U^{\rm res}_{\AA}}
\def\urese{U^{\rm res}_{\vep}}
\def\uresez{U^{\rm res}_{\vep,\ZZ}}
\def\util{\widetilde\uq}
\def\uup{U^{\geq}}
\def\ulow{U^{\leq}}
\def\bup{B^{\geq}}
\def\blow{\ovl B^{\leq}}
\def\vep{\varepsilon}
\def\vp{\varphi}
\def\vpi{\varphi^{-1}}
\def\VV{{\mathcal V}}
\def\xii{\xi^{(i)}}
\def\Xiioi{\Xi_{\io}^{(i)}}
\def\W1{W(\varpi_1)}
\def\WW{{\mathcal W}}
\def\wt{{\rm wt}}
\def\wtil{\widetilde}
\def\what{\widehat}
\def\wpi{\widehat\pi_{\lm}}
\def\ZZ{\mathbb Z}
\def\RR{\mathbb R}

\def\m@th{\mathsurround=0pt}
\def\fsquare(#1,#2){
\hbox{\vrule$\hskip-0.4pt\vcenter to #1{\normalbaselines\m@th
\hrule\vfil\hbox to #1{\hfill$\scriptstyle #2$\hfill}\vfil\hrule}$\hskip-0.4pt
\vrule}}

\theoremstyle{definition}
\newtheorem{df}{Definition}[section]
\newtheorem{thm}[df]{Theorem}
\newtheorem{pro}[df]{Proposition}
\newtheorem{lem}[df]{Lemma}
\newtheorem{ex}[df]{Example}
\newtheorem{cor}[df]{Corollary}
\newtheorem{conj}[df]{Conjecture}

\newcommand{\cmt}{\marginpar}
\newcommand{\seteq}{\mathbin{:=}}
\newcommand{\cl}{\colon}
\newcommand{\be}{\begin{enumerate}}
\newcommand{\ee}{\end{enumerate}}
\newcommand{\bnum}{\be[{\rm (i)}]}
\newcommand{\enum}{\ee}
\newcommand{\ro}{{\rm(}}
\newcommand{\rf}{{\rm)}}
\newcommand{\set}[2]{\left\{#1\,\vert\,#2\right\}}
\newcommand{\sbigoplus}{{\mbox{\small{$\bigoplus$}}}}
\newcommand{\ba}{\begin{array}}
\newcommand{\ea}{\end{array}}
\newcommand{\on}{\operatorname}
\newcommand{\eq}{\begin{eqnarray}}
\newcommand{\eneq}{\end{eqnarray}}
\newcommand{\hs}{\hspace*}

\title[Ultra-discretization 
of the $\TY(D,3,4)$-Geometric Crystals]
{Ultra-discretization 
of the $\TY(D,3,4)$-Geometric Crystals to the
$\TY(G,1,2)$-Perfect Crystals}

\author{Mana I\textsc{garashi}}
\address{Department of Mathematics, 
Sophia University, Kioicho 7-1, Chiyoda-ku, Tokyo 102-8554,
Japan}
\email{mana-i@hoffman.cc.sophia.ac.jp}
\author{Kailash C. M\textsc{isra}}
\address{Department of Mathematics,
North Carolina State University, Raleigh, NC 27695-8205, USA}
\email{misra@math.ncsu.edu}
\author{Toshiki N\textsc{akashima}}
\address{Department of Mathematics, 
Sophia University, Kioicho 7-1, Chiyoda-ku, Tokyo 102-8554,
Japan}
\email{toshiki@mm.sophia.ac.jp}
\thanks{KCM: supported in part by NSA Grant H98230-08-1-0080 and TN: supported in part by JSPS Grants 
in Aid for Scientific Research $\sharp 19540050$.}

%    General info
\subjclass{Primary 17B37; 17B67; Secondary 22E65; 14M15}
\date{}

%\dedicatory{}

\keywords{geometric crystal, perfect crystal, 
ultra-discretization. }

\begin{abstract}
Let $\ge$ be an affine Lie algebra and $\ge^L$ be its Langlands dual. It is conjectured in \cite{KNO} that  $\ge$ has a positive geometric crystal whose ultra-discretization is isomorphic to the limit of certain coherent family of perfect crystals for $\ge^L$. We prove that the ultra-discretization of the positive geometric crystal for $\ge = \TY(D,3,4)$ given in \cite{IN} is isomorphic to the limit of the coherent family of perfect crystals for $\ge^L = \TY(G,1,2)$ constructed in \cite{MMO}. 

\end{abstract}

\maketitle
%%%%%%% Sect 1. %%%%%%%%%%%%%%
\renewcommand{\thesection}{\arabic{section}}
\section{Introduction}
\setcounter{equation}{0}
\renewcommand{\theequation}{\thesection.\arabic{equation}}

Let $A= (a_{ij})_{i,j \in I}, I = \{0, 1, \cdots , n\}$ be an affine
Cartan matrix and $(A, \{\al_i\}_{i \in I}, $\\
$\{\al^\vee_i\}_{\i\in I})$ be a given Cartan datum. 
Let $\ge = \ge(A)$ denote the associated affine Lie algebra \cite{Kac} 
and $U_q(\ge)$ denote the corresponding quantum affine algebra. 
Let $P= \Z \Lambda_0 \oplus \Z \Lambda_1\oplus 
\cdots \oplus \Z \Lambda_n \oplus \Z\delta$ and  
$P^\vee = \Z \al^\vee_0 \oplus \Z \al^\vee_1 \oplus \cdots \oplus 
\Z \al^\vee_n \oplus \Z d $ denote the affine weight lattice 
and the dual affine weight lattice respectively. 
For a dominant weight  $\lambda \in P^+ = \{\mu \in P \mid \mu (h_i) 
\geq 0 \quad  {\rm for \ \  all} \quad i \in I \}$ of level 
$l = \lambda (c)$ ($c =$ canonical central element), 
Kashiwara defined the crystal base $(L(\lambda), B(\lambda))$
\cite{Kas1} 
for the integrable highest weight $U_q(\ge)$-module $V(\lambda)$. 
The crystal $B(\lambda)$ is the $q= 0$ limit of the canonical basis 
\cite{Lu} or the global crystal basis \cite{Kas2}. 
It has many interesting combinatorial properties. 
To give explicit realization of the crystal $B(\lambda)$, 
the notion of affine crystal and perfect crystal has been introduced 
in \cite{KMN1}. In particular, it is shown in \cite{KMN1} that 
the affine crystal $B(\lambda)$ for the level $l \in \Z_{>0}$ 
integrable highest weight $U_q(\ge)$-module $V(\lambda)$ can be 
realized as the semi-infinite tensor product $\cdots \otimes B_l \otimes
B_l \otimes B_l$, 
where $B_l$ is a perfect crystal of level $l$. 
This is known as the path realization. 
Subsequently it is noticed in \cite{KKM} that one needs 
a coherent family of perfect crystals $\{B_l\}_{l \geq 1}$ 
in order to give a path realization of the Verma module $M(\lambda)$ 
( or $U_q^-(\ge)$). In particular, 
the crystal $B(\infty)$ of $U_q^-(\ge)$ can be realized as the
semi-infinite tensor product $\cdots \otimes B_{\infty} \otimes 
B_{\infty} \otimes B_{\infty}$ where $B_{\infty}$ is the limit of 
the coherent family of perfect crystals $\{B_l\}_{l \geq 1}$ (see
\cite{KKM}). 
At least one coherent family $\{B_l\}_{l \geq 1}$ of perfect crystals 
and its limit is known for $\ge = A_n^{(1)}, 
B_n^{(1)}, C_n^{(1)}, D_n^{(1)}, A_{2n-1}^{(2)}, A_{2n}^{(2)}, 
D_{n+1}^{(2)}, D_4^{(3)}, G_2^{(1)} $ 
(see \cite{KMN2}, \cite{KKM}, \cite{Y}, \cite{KMOY}, \cite{MMO} ).

A perfect crystal is indeed a crystal for certain finite dimensional 
module called Kirillov-Reshetikhin module (KR-module for short) 
of the quantum affine algebra $U_q(\ge)$ (\cite{KR}, \cite{HKOTY},
\cite{HKOTT}). 
The KR-modules are parametrized by two integers $(i, l)$, where 
$i \in I \setminus \{0\}$ and $l$ any positive integer. 
Let $\{\varpi_i\}_{i\in I\setminus\{0\}}$ be the set of level $0$ 
fundamental weights \cite{K0} . Hatayama et al  
(\cite{HKOTY}, \cite{HKOTT}) conjectured that any KR-module
$W(l\varpi_i)$ 
admit a crystal base $B^{i,l}$ in the sense of Kashiwara 
and furthermore $B^{i,l}$ is perfect if $l$ is a multiple of 
$c_i^\vee\seteq \mathrm{max }(1,\frac{2}{(\al_i,\al_i)})$. 
This conjecture has been proved recently for quantum affine algebras 
$U_q(\ge)$ of classical types (\cite{OS}, \cite{FOS1}, \cite{FOS2}). 
When $\{B^{i,l}\}_{l\geq 1}$ is a coherent family of perfect crystals 
we denote its limit by $B_\infty (\varpi_i)$ (or just $B_\infty$ if there is no confusion).

On the other hand the notion of geometric crystal is introduced in
\cite{BK} 
as a geometric analog to Kashiwara's crystal (or algebraic crystal)
\cite{Kas1}. 
In fact, geometric crystal is defined in \cite{BK} for reductive
algebraic groups 
and is extended to general Kac-Moody groups in \cite{N}. 
For a given Cartan datum $(A, \{\alpha_i\}_{i \in I},
\{\al^\vee_i\}_{\i\in I} )$, 
the geometric crystal is defined as a quadruple 
$\cV(\ge)=(X, \{e_i\}_{i \in I}, \{\gamma_i\}_{i \in I}, 
\{\vep_i\}_{i\in I})$, 
where $X$ is an algebraic variety,  $e_i:\bbC^\times\times
X\longrightarrow X$ 
are rational $\bbC^\times$-actions and  
$\gamma_i,\vep_i:X\longrightarrow 
\bbC$ $(i\in I)$ are rational functions satisfying certain conditions  
( see Definition \ref{def-gc}). Geometric crystals have many properties 
similar to algebraic crystals. For instance, the product of two 
geometric crystals admits the structure of a geometric crystal if they
are induced from unipotent crystals (see \cite{BK}). A geometric 
crystal is said to be a positive geometric crystal 
if it admits a positive structure (see Definition 2.5).
A remarkable relation between positive geometric crystals 
and algebraic crystals is the ultra-discretization functor $\mathcal
{UD}$ 
between them (see Section 2.4). Applying this functor, positive rational 
functions are transfered to piecewise linear 
functions by the simple correspondence:
$$
x \times y \longmapsto x+y, \qquad \frac{x}{y} \longmapsto x - y, 
\qquad x + y \longmapsto {\rm max}\{x, y\}.
$$

Let $G$ denote the affine Kac-Moody group associated with the affine Lie
algebra $\ge$. 
Let $B^\pm$ be fixed Borel subgroups and $T$ the  maximal torus of $G$
such that $B^+\cap B^-=T$. Set $y_i(c)\seteq\exp(cf_i)$,
and let $\al_i^\vee(c)\in T$ 
be the image of $c\in\bbC^\times$
by the group morphism $\bbC^\times\to T$ induced by
the simple coroot $\alpha_i^\vee$.
We set $Y_i(c)\seteq y_i(c^{-1})\,\al_i^\vee(c)=\al_i^\vee(c)\,y_i(c)$.
Let $W$ (resp.~$\wtil W$) be the Weyl group 
(resp.~the extended Weyl group) associated with $\ge$.
The Schubert cell $X_w\seteq BwB/B$ $(w=s_{i_1}\cd s_{i_k}\in W)$ 
is birationally isomorphic to the variety
\[
 B^-_\io\seteq\set{Y_{i_1}(x_1)\cd Y_{i_k}(x_k)}%
{x_1,\cd,x_k\in \bbC^\times}\subset B^-,
\]
and $X_w$ has a natural geometric crystal structure, where 
$\io=i_1,\cd,i_k$ is a reduced word for $w$.
(\cite{BK}, \cite{N}).

Let $W(\varpi_i)$ be the KR-module 
(also called the fundamental representation)
of $U_q(\ge)$ with $\varpi_i$ as an extremal weight (\cite{K0}).
Let us denote its specialization at $q=1$
by the same notation $W(\varpi_i)$. 
It is a finite-dimensional $\ge$-module (not necessarily irreducible).
Let $\bbP(\varpi_i)$ be the projective space
$(W(\varpi_i)\setminus\{0\})/\bbC^\times$.
For any $i \in I$ the translation $t(c^\vee_i\varpi_i)$ 
belongs to $\widetilde W$ (see \cite{KNO}).
For a subset $J$ of $I$, let us denote by
$\ge_J$ the subalgebra of $\ge$ generated by
$\{e_i,f_i\}_{i\in J}$.
For an integral weight $\mu$, define 
$I(\mu)\seteq\set{j\in I}{\lan \al^\vee_j,\mu\ran\geq0}$. We recall 
the following conjecture stated in \cite{KNO}.

\begin{conj}[\cite{KNO}]
For any $i\in I \setminus \{0\}$
there exist a unique variety $X$ endowed with
a positive $\ge$-geometric crystal structure and a
rational mapping $\pi\cl X\longrightarrow
\bbP(\varpi_i)$ satisfying the following property:
\begin{enumerate}
\item
for an arbitrary extremal vector $u\in W(\varpi_i)_\mu$,
writing the translation
$t(c_i^\vee\mu)$ as $\tau w\in 
\wtil W$ with a Dynkin diagram automorphism $\tau$
and $w=s_{i_1}\cd s_{i_k}$,
%\ro see {\rm \ref{shift}}\rf, 
there exists a birational mapping
$\xi\cl B^-_{i_1,\cd,i_k}\longrightarrow X$
such that $\xi$ is a morphism of $\ge_{I(\mu)}$-geometric crystals
and that
the composition
$\pi\circ\xi\cl B^-_{i_1,\cd,i_k}\to \bbP(\varpi_i)$
coincides with
$Y_{i_1}(x_1)\cdots Y_{i_k}(x_k)\mapsto 
Y_{i_1}(x_1)\cd Y_{i_k}(x_k)\ovl u$,
where $\ovl u$ is the line including $u$,
\item
the ultra-discretization(see Sect.2) of $X$
is isomorphic to the crystal $B_\infty = B_\infty(\varpi_i)$
of the Langlands dual $\ge^L$.
\end{enumerate}
\end{conj}

In \cite{KNO},  it has been shown that this conjecture is true for $i=1$ and $\ge = A_n^{(1)}, 
B_n^{(1)}, C_n^{(1)}, D_n^{(1)}, A_{2n-1}^{(2)}, A_{2n}^{(2)},
D_{n+1}^{(2)}$. In \cite{N3}, 
a positive geometric crystal for $\ge = G_2^{(1)}$ and $i=1$ has been
constructed and it is shown in \cite{N4} that 
the  ultra-discretization of this positive geometric crystal is 
isomorphic to the limit of the coherent family of perfect crystals for 
$\ge ^L= D_4^{(3)}$ given in \cite{KMOY}. 

More recently, two of the authors have constructed a positive geometric
crystal for $\ge = D_4^{(3)}, i=1$ in \cite{IN}. 
In this paper we describe the structure of the crystal obtained 
by the ultra-discretization of the geometric crystal $\cV(\ge)$ 
constructed in \cite{IN} and then prove that it is isomorphic to the limit 
$B_\infty$ of the coherent family of perfect crystals for its 
Langlands dual $\ge ^L = G_2^{(1)}$constructed in \cite{MMO}. 
This proves Conjecture 4.5 in \cite{IN}. 

This paper is organized as follows. In Section 2, 
we recall necessary definitions and facts about geometric crystals. 
In Section 3, we review needed facts about affine crystals and perfect
crystals. 
We recall from \cite{MMO} the coherent family of perfect crystals 
for $\ge = G_2^{(1)}$ and its limit in Section 4. In Sections 5, 
we review the positive geometric crystal $\cV(\ge)$ for $\ge =
D_4^{(3)}$ constructed in \cite{IN}. 
In Section 6, we state and prove our main result (Theorem 7.1).

%%%%%%%%%%%%% Sect. 2 %%%%%%%%%%%%
\renewcommand{\thesection}{\arabic{section}}
\section{Geometric crystals}
\setcounter{equation}{0}
\renewcommand{\theequation}{\thesection.\arabic{equation}}

In this section, 
we review Kac-Moody groups and geometric crystals
following 
\cite{PK}, \cite{Ku2}, \cite{BK}
\subsection{Kac-Moody algebras and Kac-Moody groups}
\label{KM}
Fix a symmetrizable generalized Cartan matrix
 $A=(a_{ij})_{i,j\in I}$ with a finite index set $I$.
Let $(\tt,\{\al_i\}_{i\in I},\{\al^\vee_i\}_{i\in I})$ 
be the associated
root data, where ${\tt}$ is a vector space 
over $\bbC$ and
$\{\al_i\}_{i\in I}\subset\tt^*$ and 
$\{\al^\vee_i\}_{i\in I}\subset\tt$
are linearly independent 
satisfying $\al_j(\al^\vee_i)=a_{ij}$.

The Kac-Moody Lie algebra $\ge=\ge(A)$ associated with $A$
is the Lie algebra over $\bbC$ generated by $\tt$, the 
Chevalley generators $e_i$ and $f_i$ $(i\in I)$
with the usual defining relations (\cite{KP},\cite{PK}).
There is the root space decomposition 
$\ge=\bigoplus_{\al\in \tt^*}\ge_{\al}$.
Denote the set of roots by 
$\Delta:=\{\al\in \tt^*|\al\ne0,\,\,\ge_{\al}\ne(0)\}$.
Set $Q=\sum_i\bbZ \al_i$, $Q_+=\sum_i\bbZ_{\geq0} \al_i$,
$Q^\vee:=\sum_i\bbZ \al^\vee_i$
and $\Delta_+:=\Delta\cap Q_+$.
An element of $\Delta_+$ is called 
a {\it positive root}.
Let $P\subset \tt^*$ be a weight lattice such that 
$\bbC\ot P=\tt^*$, whose element is called a
weight.

Define simple reflections $s_i\in{\rm Aut}(\tt)$ $(i\in I)$ by
$s_i(h):=h-\al_i(h)\al^\vee_i$, which generate the Weyl group $W$.
It induces the action of $W$ on $\tt^*$ by
$s_i(\lm):=\lm-\lm(\al^\vee_i)\al_i$.
Set $\Delre:=\{w(\al_i)|w\in W,\,\,i\in I\}$, whose element 
is called a real root.

Let $\ge'$ be the derived Lie algebra 
of $\ge$ and let 
$G$ be the Kac-Moody group associated 
with $\ge'$(\cite{PK}).
Let $U_{\al}:=\exp\ge_{\al}$ $(\al\in \Delre)$
be the one-parameter subgroup of $G$.
The group $G$ is generated by $U_{\al}$ $(\al\in \Delre)$.
Let $U^{\pm}$ be the subgroup generated by $U_{\pm\al}$
($\al\in \Delre_+=\Delre\cap Q_+$), {\it i.e.,}
$U^{\pm}:=\lan U_{\pm\al}|\al\in\Del^{\rm re}_+\ran$.

For any $i\in I$, there exists a unique homomorphism;
$\phi_i:SL_2(\bbC)\rightarrow G$ such that
\[
\hspace{-2pt}\phi_i\left(
\left(
\begin{array}{cc}
c&0\\
0&c^{-1}
\end{array}
\right)\right)=c^{\al^\vee_i},\,
\phi_i\left(
\left(
\begin{array}{cc}
1&t\\
0&1
\end{array}
\right)\right)=\exp(t e_i),\,
 \phi_i\left(
\left(
\begin{array}{cc}
1&0\\
t&1
\end{array}
\right)\right)=\exp(t f_i).
\]
where $c\in\bbC^\times$ and $t\in\bbC$.
Set $\al^\vee_i(c):=c^{\al^\vee_i}$,
$x_i(t):=\exp{(t e_i)}$, $y_i(t):=\exp{(t f_i)}$, 
$G_i:=\phi_i(SL_2(\bbC))$,
$T_i:=\phi_i(\{{\rm diag}(c,c^{-1})\vert 
c\in\bbC^{\vee}\})$ 
and 
$N_i:=N_{G_i}(T_i)$. Let
$T$ (resp. $N$) be the subgroup of $G$ 
with the Lie algebra $\tt$
(resp. generated by the $N_i$'s), 
which is called a {\it maximal torus} in $G$, and let
$B^{\pm}=U^{\pm}T$ be the Borel subgroup of $G$.
We have the isomorphism
$\phi:W\mapright{\sim}N/T$ defined by $\phi(s_i)=N_iT/T$.
An element $\ovl s_i:=x_i(-1)y_i(1)x_i(-1)
=\phi_i\left(
\left(
\begin{array}{cc}
0&\pm1\\
\mp1&0
\end{array}
\right)\right)$ is in 
$N_G(T)$, which is a representative of 
$s_i\in W=N_G(T)/T$. 

%%%%%%%%%%%%%%%%%%%%%%%%%%%%%%%%%%%%%%%%%%%%%
\subsection{Geometric crystals}
Let $X$ be an ind-variety , 
{$\gamma_i:X\rightarrow \bbC$} and 
$\vep_i:X\longrightarrow \bbC$ ($i\in I$) 
rational functions on $X$, and
{$e_i:\bbC^\times \times X\longrightarrow X$}
$((c,x)\mapsto e^c_i(x))$ a
rational $\bbC^\times$-action.

\begin{df}
\label{def-gc}
A quadruple $(X,\{e_i\}_{i\in I},\{\gamma_i,\}_{i\in I},
\{\vep_i\}_{i\in I})$ is a 
$G$ (or $\ge$)-\\{\it geometric} {\it crystal} 
if
\begin{enumerate}
\item
$\{1\}\times X\subset dom(e_i)$ 
for any $i\in I$.
\item
$\gamma_j(e^c_i(x))=c^{a_{ij}}\gamma_j(x)$.
\item $e_i$'s satisfy the following relations.
\[
 \begin{array}{lll}
&\hspace{-20pt}e^{c_1}_{i}e^{c_2}_{j}
=e^{c_2}_{j}e^{c_1}_{i}&
{\rm if }\,\,a_{ij}=a_{ji}=0,\\
&\hspace{-20pt} e^{c_1}_{i}e^{c_1c_2}_{j}e^{c_2}_{i}
=e^{c_2}_{j}e^{c_1c_2}_{i}e^{c_1}_{j}&
{\rm if }\,\,a_{ij}=a_{ji}=-1,\\
&\hspace{-20pt}
e^{c_1}_{i}e^{c^2_1c_2}_{j}e^{c_1c_2}_{i}e^{c_2}_{j}
=e^{c_2}_{j}e^{c_1c_2}_{i}e^{c^2_1c_2}_{j}e^{c_1}_{i}&
{\rm if }\,\,a_{ij}=-2,\,
a_{ji}=-1,\\
&\hspace{-20pt}
e^{c_1}_{i}e^{c^3_1c_2}_{j}e^{c^2_1c_2}_{i}
e^{c^3_1c^2_2}_{j}e^{c_1c_2}_{i}e^{c_2}_{j}
=e^{c_2}_{j}e^{c_1c_2}_{i}e^{c^3_1c^2_2}_{j}e^{c^2_1c_2}_{i}
e^{c^3_1c_2}_je^{c_1}_i&
{\rm if }\,\,a_{ij}=-3,\,
a_{ji}=-1,
\end{array}
\]
\item
$\vep_i(e_i^c(x))=c^{-1}\vep_i(x)$ and $\vep_i(e_j^c(x))=\vep_i(x)$ if 
$a_{i,j}=a_{j,i}=0$.
\end{enumerate}
\end{df}

The condition (iv) is slightly modified from the one in 
\cite{IN,N3,N4}.

Let $W$ be the  Weyl group associated with $\ge$. 
Define $R(w)$ for $w\in W$ by
\[
 R(w):=\{(i_1,i_2,\cd,i_l)\in I^l|w=s_{i_1}s_{i_2}\cd s_{i_l}\},
\]
where $l$ is the length of $w$.
Then $R(w)$ is the set of reduced words of $w$.
For a word ${\bf i}=(i_1,\cd,i_l)\in R(w)$ 
$(w\in W)$, set 
$\al^{(j)}:=s_{i_l}\cd s_{i_{j+1}}(\al_{i_j})$ 
$(1\leq j\leq l)$ and 
\begin{eqnarray*}
e_{\bf i}:&T\times X\rightarrow &X\\
&(t,x)\mapsto &e_{\bf i}^t(x):=e_{i_1}^{\al^{(1)}(t)}
e_{i_2}^{\al^{(2)}(t)}\cd e_{i_l}^{\al^{(l)}(t)}(x).
\label{tx}
\end{eqnarray*}
Note that the condition (iii) above is 
equivalent to the following:
{$e_{\bf i}=e_{\bf i'}$}
for any 
$w\in W$, ${\bf i}$.
${\bf i'}\in R(w)$.

%%%%%%%%%%%%%%%%%%%%%%%%%%%%%%%%%%%%
\subsection{Geometric crystal on Schubert cell}
\label{schubert}

Let $w\in W$ be a Weyl group element and take a 
reduced expression $w=s_{i_1}\cd s_{i_l}$. 
Let $X:=G/B$ be the flag
variety, which is an ind-variety 
and $X_w\subset X$ the
Schubert cell associated with $w$, which has 
a natural geometric crystal structure
(\cite{BK},\cite{N}).
For ${\bf i}:=(i_1,\cd,i_k)$, set 
\begin{equation}
B_{\bf i}^-
:=\{Y_{\bf i}(c_1,\cd,c_k)
:=Y_{i_1}(c_1)\cd Y_{i_l}(c_k)
\,\vert\, c_1\cd,c_k\in\bbC^\times\}\subset B^-,
%(Y_i(c):=y_i(\frac{1}{c})\al_i^\vee(c)),
\label{bw1}
\end{equation}
where $Y_i(c):=y_i(\frac{1}{c})\al^\vee_i(c)$.
This has a geometric crystal structure(\cite{N})
isomorphic to $X_w$. 
The explicit forms of the action $e^c_i$, the rational 
function $\vep_i$  and $\gamma_i$ on 
$B_{\bf i}^-$ are given by
\begin{eqnarray}
&& e_i^c(Y_{\bf i}(c_1,\cd,c_k))
=Y_{\bf i}({\mathcal C}_1,\cd,{\mathcal C}_k)),\nn \\
&&\text{where}\nn\\
&&{\mathcal C}_j:=
c_j\cdot \frac{\displaystyle \sum_{1\leq m\leq j,i_m=i}
 \frac{c}
{c_1^{a_{i_1,i}}\cd c_{m-1}^{a_{i_{m-1},i}}c_m}
+\sum_{j< m\leq k,i_m=i} \frac{1}
{c_1^{a_{i_1,i}}\cd c_{m-1}^{a_{i_{m-1},i}}c_m}}
{\displaystyle\sum_{1\leq m<j,i_m=i} 
 \frac{c}
{c_1^{a_{i_1,i}}\cd c_{m-1}^{a_{i_{m-1},i}}c_m}+
\mathop\sum_{j\leq m\leq k,i_m=i}  \frac{1}
{c_1^{a_{i_1,i}}\cd c_{m-1}^{a_{i_{m-1},i}}c_m}},
\label{eici}\\
&& \vep_i(Y_{\bf i}(c_1,\cd,c_k))=
\sum_{1\leq m\leq k,i_m=i} \frac{1}
{c_1^{a_{i_1,i}}\cd c_{m-1}^{a_{i_{m-1},i}}c_m},
\label{vep-i}\\
&&\gamma_i(Y_{\bf i}(c_1,\cd,c_k))
=c_1^{a_{i_1,i}}\cd c_k^{a_{i_k,i}}.
\label{gamma-i}
\end{eqnarray}

%%%%%%%%%%%%%%%%%%%%%%%%%%%%%%%%%%
\subsection{Positive structure,\,\,
Ultra-discretizations \,\, and \,\,Tropicalizations}
\label{positive-str}

Let us recall the notions of 
positive structure, ultra-discretization and tropicalization.

The setting below is same as in \cite{KNO}.
Let $T=(\bbC^\times)^l$ be an algebraic torus over $\bbC$ and 
$X^*(T):={\rm Hom}(T,\bbC^\times)\cong \ZZ^l$ 
(resp. $X_*(T):={\rm Hom}(\bbC^\times,T)\cong \ZZ^l$) 
be the lattice of characters
(resp. co-characters)
of $T$. 
Set $R:=\bbC(c)$ and define
$$
\begin{array}{cccc}
v:&R\setminus\{0\}&\longrightarrow &\ZZ\\
&f(c)&\mapsto
&{\rm deg}(f(c)),
\end{array}
$$
%For any $\phi\in L(T)$, 
%set 
%${\rm deg}_T(\phi):=v\circ\phi|_{X^*(T)}$. Since 
where $\rm deg$ is the degree of poles at $c=\ify$. 
Here note that for $f_1,f_2\in R\setminus\{0\}$, we have
\begin{equation}
v(f_1 f_2)=v(f_1)+v(f_2),\q
v\left(\frac{f_1}{f_2}\right)=v(f_1)-v(f_2)
\label{ff=f+f}
\end{equation}
A non-zero rational function on
an algebraic torus $T$ is called {\em positive} if
it can be written as $g/h$ where
$g$ and $h$ are a positive linear combination of
characters of $T$.
\begin{df}
Let 
$f\cl T\rightarrow T'$ be 
a rational morphism between
two algebraic tori $T$ and 
$T'$.
We say that $f$ is {\em positive},
if $\eta\circ f$ is positive
for any character $\eta\cl T'\to \C$.
%for each $f_i$, there exist 
%polynomials
%$g_i(x_1,\cd,x_m)$, $h_i(x_1,\cd,x_m)$
%with positive coefficients such that $f_i=g_i/h_i
%\ne0$.
\end{df}
Denote by ${\rm Mor}^+(T,T')$ the set of 
positive rational morphisms from $T$ to $T'$.

\begin{lem}[\cite{BK}]
\label{TTT}
For any $f\in {\rm Mor}^+(T_1,T_2)$             
and $g\in {\rm Mor}^+(T_2,T_3)$, 
the composition $g\circ f$
is well-defined and belongs to ${\rm Mor}^+(T_1,T_3)$.
\end{lem}

By Lemma \ref{TTT}, we can define a category ${\mathcal T}_+$
whose objects are algebraic tori over $\bbC$ and arrows
are positive rational morphisms.

Let $f\cl T\rightarrow T'$ be a 
positive rational morphism
of algebraic tori $T$ and 
$T'$.
We define a map $\what f\cl X_*(T)\rightarrow X_*(T')$ by 
\[
\langle\eta,\what f(\xi)\rangle
=v(\eta\circ f\circ \xi),
\]
where $\eta\in X^*(T')$ and $\xi\in X_*(T)$.
\begin{lem}[\cite{BK}]
For any algebraic tori $T_1$, $T_2$, $T_3$, 
and positive rational morphisms 
$f\in {\rm Mor}^+(T_1,T_2)$, 
$g\in {\rm Mor}^+(T_2,T_3)$, we have
$\what{g\circ f}=\what g\circ\what f.$
\end{lem}
%{\sl Proof.}
%For $\mu\in X_*(T_1)$, 
%Let ${\mathcal B}$ be the category of free $\ZZ$-modules,
%whose arrows are piece-wise linear maps.
Let ${\hbox{\germ Set}}$ denote the category of sets with the morphisms being set maps.
By the above lemma, we obtain a functor: 
\[
\begin{array}{cccc}
{\mathcal UD}:&{\mathcal T}_+&\longrightarrow &{{\hbox{\germ Set}}}\\
&T&\mapsto& X_*(T)\\
&(f:T\rightarrow T')&\mapsto& 
(\what f:X_*(T)\rightarrow X_*(T')))
\end{array}
\]

%More precisely, for a geometric crystal $X$
%$(\gamma:X\rightarrow T)$ and a 
%positive structure $\theta:T'\rightarrow X$, 
%functor ${\rm Trop}_\theta:X\mapsto X_*(T')$ is called
%``tropicalization'' in \cite{BK}. But we call it a
%``ultra discretization (UD)''.
%And if some free crystal $B$ 
%is obtained by UD of a geometric crystal $X$, we call
%$X$ a {\it tropicalization} of $B$.

\begin{df}[\cite{BK}]
Let $\chi=(X,\{e_i\}_{i\in I},\{{\rm wt}_i\}_{i\in I},
\{\vep_i\}_{i\in I})$ be a 
geometric crystal, $T'$ an algebraic torus
and $\theta:T'\rightarrow X$ 
a birational isomorphism.
The isomorphism $\theta$ is called 
{\it positive structure} on
$\chi$ if it satisfies
\begin{enumerate}
\item for any $i\in I$ the rational functions
$\gamma_i\circ \theta:T'\rightarrow \bbC$ and 
$\vep_i\circ \theta:T'\rightarrow \bbC$ 
are positive.
\item
For any $i\in I$, the rational morphism 
$e_{i,\theta}:\bbC^\tm \tm T'\rightarrow T'$ defined by
$e_{i,\theta}(c,t)
:=\theta^{-1}\circ e_i^c\circ \theta(t)$
is positive.
\end{enumerate}
\end{df}
Let $\theta:T\rightarrow X$ be a positive structure on 
a geometric crystal $\chi=(X,\{e_i\}_{i\in I},$
$\{{\rm wt}_i\}_{i\in I},
\{\vep_i\}_{i\in I})$.
Applying the functor ${\mathcal UD}$ 
to positive rational morphisms
$e_{i,\theta}:\bbC^\tm \tm T'\rightarrow T'$ and
$\gamma\circ \theta:T'\ra T$
(the notations are
as above), we obtain
\begin{eqnarray*}
\til e_i&:=&{\mathcal UD}(e_{i,\theta}):
\ZZ\tm X_*(T) \rightarrow X_*(T)\\
{\rm wt}_i&:=&{\mathcal UD}(\gamma_i\circ\theta):
X_*(T')\rightarrow \bbZ,\\
\vep_i&:=&{\mathcal UD}(\vep_i\circ\theta):
X_*(T')\rightarrow \bbZ.
\end{eqnarray*}
Now, for given positive structure $\theta:T'\rightarrow X$
on a geometric crystal 
$\chi=(X,\{e_i\}_{i\in I},$
$\{{\rm wt}_i\}_{i\in I},
\{\vep_i\}_{i\in I})$, we associate 
the quadruple $(X_*(T'),\{\til e_i\}_{i\in I},
\{{\rm wt}_i\}_{i\in I},\{\vep_i\}_{i\in I})$
with a free pre-crystal structure (see \cite[2.2]{BK}) 
and denote it by ${\mathcal UD}_{\theta,T'}(\chi)$.
We have the following theorem:

\begin{thm}[\cite{BK}\cite{N}]
For any geometric crystal 
$\chi=(X,\{e_i\}_{i\in I},\{\gamma_i\}_{i\in I},$
$\{\vep_i\}_{i\in I})$ and positive structure
$\theta:T'\rightarrow X$, the associated pre-crystal
${\mathcal UD}_{\theta,T'}(\chi)=$\\
$(X_*(T'),\{e_i\}_{i\in I},\{{\rm wt}_i\}_{i\in I},
\{\vep_i\}_{i\in I})$ 
is a crystal {\rm (see \cite[2.2]{BK})}
\end{thm}

Now, let ${\mathcal GC}^+$ be a category whose 
object is a triplet
$(\chi,T',\theta)$ where 
$\chi=(X,\{e_i\},\{\gamma_i\},\{\vep_i\})$ 
is a geometric crystal and $\theta:T'\rightarrow X$ 
is a positive structure on $\chi$, and morphism
$f:(\chi_1,T'_1,\theta_1)\longrightarrow 
(\chi_2,T'_2,\theta_2)$ is given by a morphism 
$\vp:X_1\longrightarrow X_2$  
($\chi_i=(X_i,\cd)$) such that 
\[
f:=\theta_2^{-1}\circ\vp\circ\theta_1:T'_1\longrightarrow T'_2,
\]
is a positive rational morphism. Let ${\mathcal CR}$
be a category of crystals. 
Then by the theorem above, we have
\begin{cor}
\label{cor-posi}
The map $ \mathcal UD = \mathcal UD_{\theta,T'}$ defined above is a functor
\begin{eqnarray*}
 {\mathcal UD}&:&{\mathcal GC}^+\longrightarrow {\mathcal CR},\\
&&(\chi,T',\theta)\mapsto X_*(T'),\\
&&(f:(\chi_1,T'_1,\theta_1)\rightarrow 
(\chi_2,T'_2,\theta_2))\mapsto
(\what f:X_*(T'_1)\rightarrow X_*(T'_2)).
\end{eqnarray*}

\end{cor}
We call the functor $\mathcal UD$
{\it ``ultra-discretization''} as \cite{N},\cite{N2}
instead of ``tropicalization'' as in \cite{BK}.
And 
for a crystal $B$, if there
exists a geometric crystal $\chi$ and a positive 
structure $\theta:T'\rightarrow X$ on $\chi$ such that 
${\mathcal UD}(\chi,T',\theta)\cong B$ as crystals, 
we call an object $(\chi,T',\theta)$ in ${\mathcal GC}^+$
a {\it tropicalization} of $B$, where 
it is not known that this correspondence is a functor.

%%%%%%%%%%% section %%%%%%%%%%%%%%%
\renewcommand{\thesection}{\arabic{section}}
\section{Limit of perfect crystals}
\label{limit}
\setcounter{equation}{0}
\renewcommand{\theequation}{\thesection.\arabic{equation}}
%%%%%%%%%%%%%%%%%%%%%%%%%%%%%%%%%%%%%%%%%%%%%%%%%%%%%%%%
We review limit of perfect crystals following \cite{KKM}.
(See also \cite{KMN1},\cite{KMN2}).

%%%%%%%%%%%%%%%%%%%%%%%%%%%%%%%
\subsection{Crystals}

First we review the theory of crystals,
which is the notion obtained by
abstracting the combinatorial 
properties of crystal bases.
Let 
$(A,\{\al_i\}_{i\in I},\{\al^\vee_i\}_{\i\in I})$ be a 
Cartan data. 
\begin{df}
A {\it crystal} $B$ is a set endowed with the following maps:
\begin{eqnarray*}
&& {\rm wt}:B\lar P,\\
&&\vep_i:B\lar\ZZ\sqcup\{-\infty\},\q
  \vp_i:B\lar\ZZ\sqcup\{-\infty\} \q{\hbox{for}}\q i\in I,\\
&&\eit:B\sqcup\{0\}\lar B\sqcup\{0\},
\q\fit:B\sqcup\{0\}\lar B\sqcup\{0\}\q{\hbox{for}}\q i\in I,\\
&&\eit(0)=\fit(0)=0.
\end{eqnarray*}
Those maps satisfy the following axioms: for
 all $b,b_1,b_2 \in B$, we have
\begin{eqnarray*}
&&\vp_i(b)=\vep_i(b)+\lan \al^\vee_i,{\rm wt}
(b)\ran,\\
&&\wt(\eit b)=\wt(b)+\al_i{\hbox{ if  }}\eit b\in B,\\
&&\wt(\fit b)=\wt(b)-\al_i{\hbox{ if  }}\fit b\in B,\\
&&\eit b_2=b_1 \Longleftrightarrow \fit b_1=b_2\,\,(\,b_1,b_2 \in B),\\
&&\vep_i(b)=-\ify
   \Longrightarrow \eit b=\fit b=0.
\end{eqnarray*}
\end{df}
The following tensor product structure 
is one of the most crucial properties of crystals.
\begin{thm}
\label{tensor}
Let $B_1$ and $B_2$ be crystals.
Set
$B_1\ot B_2:=
\{b_1\otimes b_2;\;b_j\in B_j\;(j=1,2)\}$. Then we have 
\begin{enumerate}
\item $B_1\ot B_2$ is a crystal.
\item
For $b_1\in B_1$ and $b_2\in B_2$, we have
$$
\tilde f_i(b_1\otimes b_2)=
\left\{\begin{array}{ll}\tilde f_ib_1\otimes b_2&
{\rm if}\;\varphi_i(b_1)>\vep_i(b_2),\\
b_1\otimes\tilde f_ib_2&{\rm if}\;
\varphi_i(b_1)\leq\vep_i(b_2).
\end{array}\right.
$$
$$
\tilde e_i(b_1\otimes b_2)=\left\{\begin{array}{ll}
b_1\otimes \tilde e_ib_2&
{\rm if}\;\varphi_i(b_1)<\vep_i(b_2),\\
\tilde e_ib_1\otimes b_2
&{\rm if}\;\varphi_i(b_1)\geq\vep_i(b_2),
\end{array}\right.
$$
\end{enumerate}
\end{thm}

\begin{df}
Let $B_1$ and $B_2$ be crystals. A {\it strict morphism} of crystals
$\psi:B_1\lar B_2$ is a map
$\psi:B_1\sqcup\{0\} \lar B_2\sqcup\{0\}$
satisfying: 
$\psi(0)=0$, $\psi(B_1)\subset B_2$,
$\psi$ commutes with all $\eit$ and $\fit$
and
\[
\hspace{-30pt}\wt(\psi(b))=\wt(b),\q \vep_i(\psi(b))=\vep_i(b),\q
  \vp_i(\psi(b))=\vp_i(b)
\text{ for any }b\in B_1.
\]
In particular, 
a bijective strict morphism is called an 
{\it isomorphism of crystals}. 
\end{df}

\begin{ex}
\label{ex-tlm}
If $(L,B)$ is a crystal base, then $B$ is a crystal.
Hence, for the crystal base $(L(\ify),B(\ify))$
of the nilpotent subalgebra $\uqm$ of 
the quantum algebra $\uq$, 
$B(\ify)$ is a crystal. 
\end{ex}
\begin{ex}
\label{tlm}
For $\lm\in P$, set $T_\lm:=\{t_\lm\}$. We define a crystal
structure on $T_\lm$ by 
\[
 \eit(t_\lm)=\fit(t_\lm)=0,\q\vep_i(t_\lm)=
\vp_i(t_\lm)=-\ify,\q \wt(t_\lm)=\lm.
\]
\end{ex}
\begin{df}
For a crystal $B$, a colored oriented graph
structure is associated with $B$ by 
\[
 b_1\mapright{i}b_2\Longleftrightarrow 
\fit b_1=b_2.
\]
We call this graph a {\it crystal graph}
of $B$.
\end{df}

%%%%%%%%%%%%%%%%%%%%%%%%%%%%%%
\subsection{Affine weights}
\label{aff-wt}

Let $\ge$ be an affine Lie algebra. 
The sets $\mathfrak t$, 
$\{\al_i\}_{i\in I}$ 
and $\{\al^\vee_i\}_{i\in I}$ be as in \ref{KM}. 
We take ${\rm dim}\mathfrak t=\sharp I+1$.
Let $\del\in Q_+$ be the unique element 
satisfying $\{\lm\in Q|\lan \al^\vee_i,\lm\ran=0
\text{ for any }i\in I\}=\bbZ\del$
and ${\bf c}\in \ge$ be the canonical central element
satisfying $\{h\in Q^\vee|\lan h,\al_i\ran=0
\text{ for any }i\in I\}=\bbZ c$.
We write (\cite[6.1]{Kac})
\[
{\bf c}=\sum_i a_i^\vee \al^\vee_i,\qq
\del=\sum_i a_i\al_i.
\]
Let $(\q,\q)$ be the non-degenerate
$W$-invariant symmetric bilinear form on $\mathfrak t^*$
normalized by $(\del,\lm)=\lan {\bf c},\lm\ran$
for $\lm\in\frak t^*$.
Let us set $\tt^*_{\rm cl}:=\tt^*/\bbC\del$ and let
${\rm cl}:\tt^*\longrightarrow \tt^*_{\rm cl}$
be the canonical projection. 
Here we have 
$\tt^*_{\rm cl}\cong \oplus_i(\bbC \al^\vee_i)^*$.
Set $\tt^*_0:=\{\lm\in\tt^*|\lan {\bf c},\lm\ran=0\}$,
$(\tt^*_{\rm cl})_0:={\rm cl}(\tt^*_0)$. 
Since $(\del,\del)=0$, we have a positive-definite
symmetric form on $\tt^*_{\rm cl}$ 
induced by the one on 
$\tt^*$. 
Let $\Lm_i\in \tt^*_{\rm cl}$ $(i\in I)$ be a classical 
weight such that $\lan \al^\vee_i,\Lm_j\ran=\del_{i,j}$, which 
is called a fundamental weight.
We choose 
$P$ so that $P_{\rm cl}:={\rm cl}(P)$ coincides with 
$\oplus_{i\in I}\bbZ\Lm_i$ and 
we call $P_{\rm cl}$ a 
{\it classical weight lattice}.

%%%%%%%%%%%% subsection %%%%%%%%%%%
\subsection{Definitions of perfect crystal and its limit}
\label{def-perfect}

Let $\ge$ be an affine Lie algebra, $P_{cl}$ be
a classical weight lattice as above and set 
$(P_{cl})^+_l:=\{\lm\in P_{cl}|
\lan c,\lm\ran=l,\,\,\lan \al^\vee_i,\lm\ran\geq0\}$ 
$(l\in\ZZ_{>0})$.
\begin{df}
\label{perfect-def}
A crystal $B$ is a {\it perfect crystal} of level $l$ if 
\begin{enumerate}
\item
$B\ot B$ is connected as a crystal graph.
\item
There exists $\lm_0\in P_{\rm cl}$ such that 
\[
 \wt(B)\subset \lm_0+\sum_{i\ne0}\ZZ_{\leq0}
{\rm cl}(\al_i),\qq
\sharp B_{\lm_0}=1
\]
\item There exists a finite-dimensional 
$U'_q(\ge)$-module $V$ with a
crystal pseudo-base $B_{ps}$ 
such that $B\cong B_{ps}/{\pm1}$
\item
The maps 
$\vep,\vp:B^{min}:=\{b\in B|\lan c,\vep(b)\ran=l\}
\mapright{}(P_{\rm cl}^+)_l$ are bijective, where 
$\vep(b):=\sum_i\vep_i(b)\Lm_i$ and 
$\vp(b):=\sum_i\vp_i(b)\Lm_i$.
\end{enumerate}
\end{df}

Let $\{B_l\}_{l\geq1}$ be a family of 
perfect crystals of level $l$ and set 
$J:=\{(l,b)|l>0,\,b\in B^{min}_l\}$.
\begin{df}
\label{def-limit}
A crystal $B_\ify$ with an element $b_\ify$ is called a
{\it limit of $\{B_l\}_{l\geq1}$}
if 
\begin{enumerate}
\item
$\wt(b_\ify)=\vep(b_\ify)=\vp(b_\ify)=0$.
\item
For any $(l,b)\in J$, there exists an
embedding of 
crystals:
\begin{eqnarray*}
 f_{(l,b)}:&
T_{\vep(b)}\ot B_l\ot T_{-\vp(b)}\hookrightarrow
B_\ify\\
&t_{\vep(b)}\ot b\ot t_{-\vp(b)}\mapsto b_\ify
\end{eqnarray*}
\item
$B_\ify=\bigcup_{(l,b)\in J} {\rm Im}f_{(l,b)}$.
\end{enumerate}
\end{df}
\noindent
As for the crystal $T_\lm$, see Example \ref{tlm}.
If a limit exists for a family $\{B_l\}$, 
we say that $\{B_l\}$
is a {\it coherent family} of perfect crystals.

The following is one of the most 
important properties of limit of perfect crystals.
\begin{pro}
Let $B(\ify)$ be the crystal as in 
Example \ref{ex-tlm}. Then we have
the following isomorphism of crystals:
\[
B(\ify)\ot B_\ify\mapright{\sim}B(\ify).
\]
\end{pro}

%%%%%%%%%%%%%
\section{Perfect Crystals of type $\TY(G,1,2)$}
\label{perf}

In this section, we review the family of 
perfect crystals of type 
$\TY(G,1,2)$ and its limit(\cite{MMO}).

We fix the data for $G_2^{(1)}$.
Let 
$\{\alpha_0, \alpha_1, \alpha_2\}$, 
$\{\al^\vee_0, \al^\vee_1, \al^\vee_2\}$ and 
$\{\Lm_0, \Lm_1, \Lm_2\}$ be the set of 
simple roots, simple coroots and fundamental weights, respectively.
The Cartan matrix $A=(a_{ij})_{i,j=0,1,2}$
is given by
\[A=
\left(
\begin{array}{rrr}
2  & -1 & 0  \\
-1 & 2  & -1 \\
0  & -3 & 2
\end{array}
\right),
\]
and its Dynkin diagram is as follows.
%\[
%{}^0 \circ-{}^1 \circ \Lleftarrow \circ^2
%\]
\[\SelectTips{cm}{}
\xymatrix{
*{\bigcirc}<3pt> \ar@{-}[r]_<{0} 
%& *{\bigcirc}<3pt> \ar@3{<-}[r]_<{1}
& *{\bigcirc}<3pt> \ar@3{->}[r]_<{1}
& *{\bigcirc}<6pt>\ar@{}_<{\,\,\,\,\,\,2}
}
\]
The standard null root $\delta$ 
and the canonical central element $c$ are 
given by
\[
\delta=\alpha_0+2\alpha_1+3\alpha_2
\quad\text{and}\quad c=\al^\vee_0+2\al^\vee_1+\al^\vee_2, 
\]
where 
$\al_0=2\Lm_0-\Lm_1+\del,\q
\al_1=-\Lm_0+2\Lm_1-3\Lm_2,\q
\al_2=-\Lm_1+2\Lm_2.$

For a positive integer $l$ we introduce 
$\TY(G,1,2)$-crystals $B_l$ and $B_\ify$ as 
\begin{eqnarray*}
&&B_l=\left\{
b=(b_1,b_2,b_3,{\bar b}_3,{\bar b}_2,{\bar b}_1)
\in(\Zn/3)^6
\left\vert
\begin{array}{l}
3b_3\equiv3\bar{b}_3\;(\text{mod }2), \\
\sum_{i=1,2} (b_i+{\bar b}_i)+\frac{b_3+{\bar b}_3}{2}
\leq l\\
b_1, {\bar b}_1, b_2-b_3, {\bar b}_3-{\bar b}_2 \in \ZZ
\end{array}
\right.
\right\},\\
&&B_\ify=\left\{
b=(b_1,b_2,b_3,{\bar b}_3,{\bar b}_2,{\bar b}_1)
\in(\ZZ/3)^6
\left\vert
\begin{array}{l}
3b_3\equiv3\bar{b}_3\;(\text{mod }2), \\
b_1, {\bar b}_1, b_2-b_3, {\bar b}_3-{\bar b}_2 \in \ZZ
%\sum_{i=1,2} (b_i+{\bar b}_i)+\frac{b_3+{\bar b}_3}{2}
%\in\ZZ
\end{array}
\right.
\right\}. 
\end{eqnarray*}
Now we describe the explicit crystal structures
of $B_l$ and $B_\ify$. 
Indeed, most of them coincide with 
each other except
for $\vep_0$ and $\vp_0$.
In the rest of this section, we use the following 
convention: 
$(x)_+=\max(x,0)$.
For $b=(b_1,b_2,b_3,{\bar b}_3,{\bar b}_2,{\bar b}_1)$ we denote
\begin{equation} \label{def s(b)}
s(b)=b_1+b_2+\frac{b_3+{\bar b}_3}{2}+{\bar b}_2+{\bar b}_1, 
\end{equation}
and 
\begin{equation} \label{z1-4}
z_1={\bar b}_1-b_1, \quad 
z_2={\bar b}_2 -{\bar b}_3, \quad 
z_3=b_3-b_2, \quad 
z_4=({\bar b}_3-b_3)/2.
\end{equation}

Now we define conditions ($E_1$)-($E_6$) and ($F_1$)-($F_6$) as follows.

\begin{equation} \label{(F)}
%\begin{split}
\begin{cases}
&(F_1)\quad
z_1+z_2+z_3+3z_4\le0, z_1+z_2+3z_4\le0, z_1+z_2\le0, z_1\le0,
\\
&(F_2)\quad
z_1+z_2+z_3+3z_4\le0, z_2+3z_4\le0, z_2\le0, z_1> 0,
\\
&(F_3)\quad 
z_1+z_3+3z_4\le0, z_3+3z_4\le0, z_4\le0, z_2> 0, z_1+z_2> 0,
\\
&(F_4)\quad
z_1+z_2+3z_4> 0, z_2+3z_4> 0, z_4> 0, z_3\le0, z_1+z_3\le0,
\\
&(F_5)\quad
z_1+z_2+z_3+3z_4> 0, z_3+3z_4> 0, z_3> 0, z_1\le0, 
\\
&(F_6)\quad
z_1+z_2+z_3+3z_4> 0, z_1+z_3+3z_4> 0, z_1+z_3> 0, z_1> 0.
%\end{split}
\end{cases}
\end{equation}
($E_i$) ($1\le i\le 6$) is defined from ($F_i$) by replacing $>$ (resp. $\le$) with $\geq$ (resp. $<$). We also define 
\begin{equation} \label{A}
A=(0,z_1,z_1+z_2,z_1+z_2+3z_4,z_1+z_2+z_3+3z_4,2z_1+z_2+z_3+3z_4).
\end{equation}

Then for $b=(b_1,b_2,b_3,{\bar b}_3,{\bar b}_2,{\bar b}_1) \in B_l$ or $B_\ify$,
$\et{i}b, \ft{i}b, \vep_i(b), \vp_i(b), i= 0, 1, 2$ are given as follows.

\begin{align*}
\et{0}b=&
\begin{cases}
(b_1 -1,\ldots) 
& \text{if ($E_1$)}, 
\\
(\ldots,b_3 -1,{\bar b}_3 -1,\ldots,{\bar b}_1 +1) 
& \text{if ($E_2$)}, 
\\
(\ldots,b_2-\frac{2}{3},b_3-\frac{2}{3},{\bar b}_3+\frac{4}{3},{\bar b}_2+\frac{1}{3},\ldots) 
& \text{if ($E_3$) and $z_4=-\frac{1}{3}$}, 
\\ 
(\ldots,b_2-\frac{1}{3},b_3-\frac{4}{3},{\bar b}_3+\frac{2}{3},{\bar b}_2+\frac{2}{3},\ldots) 
& \text{if ($E_3$) and $z_4=-\frac{2}{3}$}, 
\\ 
(\ldots,b_3 -2,\ldots,{\bar b}_2 +1,\ldots) 
& \text{if ($E_3$) and $z_4\ne-\frac{1}{3},-\frac{2}{3}$}, 
\\ 
(\ldots,b_2 -1,\ldots,{\bar b}_3 +2,\ldots) 
& \text{if ($E_4$)}, 
\\
(b_1 -1,\ldots,b_3 +1,{\bar b}_3 +1,\ldots) 
& \text{if ($E_5$)}, 
\\
(\ldots,{\bar b}_1 +1) 
& \text{if ($E_6$)},
\end{cases}
\\
\ft{0}b=&
\begin{cases}
(b_1 +1,\ldots) 
& \text{if ($F_1$)}, 
\\ 
(\ldots,b_3 +1,{\bar b}_3 +1,\ldots,{\bar b}_1 -1) 
& \text{if ($F_2$)}, 
\\
(\ldots,b_3 +2,\ldots,{\bar b}_2 -1,\ldots) 
& \text{if ($F_3$)}, 
\\
(\ldots,b_2 +\frac{1}{3},b_3+\frac{4}{3},{\bar b}_3-\frac{2}{3},{\bar b}_2 -\frac{2}{3},\ldots) 
& \text{if ($F_4$) and $z_4=\frac{1}{3}$}, 
\\
(\ldots,b_2 +\frac{2}{3},b_3+\frac{2}{3},{\bar b}_3-\frac{4}{3},{\bar b}_2 -\frac{1}{3},\ldots) 
& \text{if ($F_4$) and $z_4=\frac{2}{3}$}, 
\\
(\ldots,b_2 +1,\ldots,{\bar b}_3 -2,\ldots) 
& \text{if ($F_4$) and $z_4\ne\frac{1}{3},\frac{2}{3}$}, 
\\
(b_1 +1,\ldots,b_3 -1,{\bar b}_3 -1,\ldots) 
& \text{if ($F_5$)}, 
\\
(\ldots,{\bar b}_1 -1)
& \text{if ($F_6$)},
\end{cases}
\end{align*}

\begin{align*}
\et{1}b=&
\begin{cases}
(\ldots,{\bar b}_2 +1,{\bar b}_1 -1) 
& \text{if ${\bar b}_2 -{\bar b}_3 \geq (b_2 -b_3)_+$}, 
\\  
(\ldots,b_3 +1,{\bar b}_3 -1,\ldots) 
& \text{if ${\bar b}_2 -{\bar b}_3 <0\leq b_3 -b_2$}, 
\\ 
(b_1 +1,b_2 -1,\ldots) 
& \text{if $({\bar b}_2 -{\bar b}_3)_+ <b_2 -b_3$},
\end{cases}
\\
\ft{1}b=&
\begin{cases}
(b_1 -1,b_2 +1,\ldots) 
& \text{if $({\bar b}_2 -{\bar b}_3)_+ \leq b_2 -b_3$}, 
\\
(\ldots,b_3 -1,{\bar b}_3 +1,\ldots) 
& \text{if ${\bar b}_2 -{\bar b}_3 \leq 0< b_3 -b_2$}, 
\\
(\ldots,{\bar b}_2 -1,{\bar b}_1 +1) 
& \text {if ${\bar b}_2 -{\bar b}_3 >(b_2 -b_3)_+$},
\end{cases}
\\
\et{2}b=&
\begin{cases}
(\ldots,{\bar b}_3 +\frac23,{\bar b}_2 -\frac13,\ldots) 
& \text{if ${\bar b}_3 \geq b_3$}, 
\\
(\ldots,b_2 +\frac13,b_3 -\frac23,\ldots) 
& \text{if ${\bar b}_3 <b_3$},
\end{cases}
\\
\ft{2}b=&
\begin{cases}
(\ldots,b_2 -\frac13,b_3 +\frac23,\ldots) 
& \text{if ${\bar b}_3 \leq b_3$}, 
\\
(\ldots,{\bar b}_3 -\frac23,{\bar b}_2 +\frac13,\ldots) 
& \text{if ${\bar b}_3 >b_3$}.
\end{cases}
\end{align*}

\begin{align*}
\vep_1(b)=&{\bar b}_1+({\bar b}_3-{\bar b}_2+(b_2-b_3)_+)_+,
\qq\vp_1(b)=
b_1+(b_3-b_2+({\bar b}_2-{\bar b}_3)_+)_+,
\\
\vep_2(b)=&{3\bar b}_2+\frac{3}{2}(b_3-{\bar b}_3)_+,
\qq\vp_2(b)=3b_2+\frac{3}{2}({\bar b}_3-b_3)_+,\\
\vep_0(b)=&
\begin{cases}
l-s(b)+\max \,A-(2z_1+z_2+z_3+3z_4)&b\in B_l,\\
-s(b)+\max \,A-(2z_1+z_2+z_3+3z_4)&b\in B_\ify.
\end{cases}\\
\vp_0(b)=&
\begin{cases}
l-s(b)+\max\, A&b\in B_l,\\
-s(b)+\max\, A&b\in B_\ify,
\end{cases}
\end{align*}

For $b\in B_l$ 
if $\et{i}b$ or $\ft{i}b$ does not belong to 
$B_l$, namely, if $b_j$ or $\bar{b}_j$
for some $j$ becomes negative or $s(b)$ exceeds $l$, 
we understand it to be $0$.

The following is one of the main results in 
\cite{MMO}:
\begin{thm}[\cite{MMO}]
\begin{enumerate}
\item
The $\TY(G,1,2)$-crystal $B_l$ 
is a perfect crystal of level $l$. 
\item
The family of the perfect crystals 
$\{B_l\}_{\l\geq1}$ forms a 
coherent family and the crystal $B_\ify$ 
is its limit with the vector 
$b_\ify=(0,0,0,0,0,0)$.
\end{enumerate}
\end{thm}
As was shown in \cite{MMO}, 
the minimal elements are given
\[
(B_l)_{\min}=
\{(\alpha,\beta,\beta,\beta,\beta,\alpha)\,|\,
\al\in\ZZ_{\geq0}, \beta\in (\ZZ_{\geq0})/3, 2\al+3\beta\leq l\}.
\] 
Let $J=\{(l,b)\,|\,l\in\ZZ_{\geq1},
b\in(B_l)_{\min}\}$ and the maps $\vep,\,\vp
:(B_l)_{\min}
\to (P^+_{\rm cl})_l$ be as in Sect.3.
Then we have 
$\wt b_\infty=0$ and $\vep_i(b_\infty)
=\vp_i(b_\infty)=0$ for $i=0,1,2$.

For $(l,b_0)\in J$, since $\vep(b_0)=\vp(b_0)$, 
one can set $\lm=\vep(b_0)=\vp(b_0)$. For 
$b=(b_1,b_2,b_3,\bar{b}_3,\bar{b}_2,\bar{b}_1)\in B_l$ we define a map
\[
f_{(l,b_0)}:\;T_{\lm}\ot B_l\ot 
B_{-\lm}\longrightarrow B_\infty
\]
by
\[
f_{(l,b_0)}(t_{\lm}\ot b\ot t_{-\lm})=b'=(\nu_1,\nu_2,\nu_3,\bar{\nu}_3,\bar{\nu}_2,\bar{\nu}_1)
\]
where
$b_0=(\al,\beta,\beta,\beta,\beta,\al)$, and
\begin{align*}
\nu_1&=b_1-\alpha, & \bar{\nu}_1&=\bar{b}_1-\alpha, \\
\nu_j&=b_j-\beta, & \bar{\nu}_j&=\bar{b}_j-\beta\;(j=2,3).
\end{align*}
Finally, we obtain 
$B_\infty=\bigcup_{(l,b)\in J}
\mbox{Im}\,f_{(l,b)}$

%%%%%%%%%%%%%% Section  %%%%%%%%%%%%%%%%
\renewcommand{\thesection}{\arabic{section}}
\section{Affine Geometric Crystal $\cV_1(\TY(D,3,4))$}
\setcounter{equation}{0}
\renewcommand{\theequation}{\thesection.\arabic{equation}}

%%%%%%%%%% 
\subsection{Fundamental representation 
$W(\varpi_1)$ for $\TY(D,3,4)$}
\label{fundamental}

Let $c=\sum_{i}a_i^\vee \al^\vee_i$ be the canonical
central element in an affine Lie algebra $\ge$
(see \cite[6.1]{Kac}), 
$\{\Lm_i|i\in I\}$ the set of fundamental 
weight as in the previous section
and $\varpi_1:=\Lm_1-a^\vee_1\Lm_0$ the
(level 0)fundamental weight.
Let $W(\varpi_1)$ be the fundamental representation 
of $\uqp$
associated with $\varpi_1$ (\cite{K0}).

By \cite[Theorem 5.17]{K0}, $W(\varpi_1)$ is a
finite-dimensional irreducible integrable 
$\uqp$-module and has a global basis
with a simple crystal. Thus, we can consider 
the specialization $q=1$ and obtain the 
finite-dimensional $\ge$-module $W(\varpi_1)$, 
which we call a fundamental representation
of $\ge$ and use the same notation as above.

We shall present the explicit form of 
$W(\varpi_1)$ for $\ge=\TY(D,3,4)$.
%%%%%%%%%%% subsection %%%%%%%%%%
\subsection{$W(\varpi_1)$ for $\TY(D,3,4)$}
The Cartan matrix $A=(a_{i,j})_{i,j=0,1,2}$ of type 
$\TY(D,3,4)$ is:
\[
 A=\begin{pmatrix}2&-1&0\\
-1&2&-3\\0&-1&2
\end{pmatrix}.
\]
%Cartan matrix 
%\[
% a_{i,j}=\begin{cases}
%2&\text{ if }i=j,\\
%-1&\text{ if }(i,j)=(0,1), (1,0), (1,2),\\
%-3&\text{ if }(i,j)=(2,1),\\
%0&\text{ otherwise}.
%\end{cases}
%\]
Then the simple roots are 
\[
 \al_0=2\Lm_0-\Lm_1+\del,\q
\al_1=-\Lm_0+2\Lm_1-\Lm_2,\q
\al_2=-3\Lm_1+2\Lm_2, 
\]
and the Dynkin diagram is:
\[\SelectTips{cm}{}
\xymatrix{
*{\bigcirc}<3pt> \ar@{-}[r]_<{0} 
& *{\bigcirc}<3pt> \ar@3{<-}[r]_<{1}
& *{\bigcirc}<6pt>\ar@{}_<{\,\,\,\,\,\,2}
}
\]

The $\TY(D,3,4)$-module $W(\varpi_1)$ is an 8-dimensional 
module with the basis,
\[
 \{v_1,v_2, v_3,v_0, {\emptyset}, v_{\ovl 3},
v_{\ovl 2},v_{\ovl 1}\}.
\]
The explicit form of $W(\varpi_1)$ is given 
in \cite{KMOY}.
\begin{eqnarray*}
&&\hspace{-30pt}{\rm wt}(v_1)=\Lm_1-2\Lm_0,\,\,
{\rm wt}(v_2)=-\Lm_0-\Lm_1+\Lm_2,\,\,
{\rm wt}(v_3)=-\Lm_0+2\Lm_1-\Lm_2,\\
&&\hspace{-30pt}{\rm wt}(v_{\ovl i})=
-{\rm wt}(v_i)\,\,(i=1,\cd,3),\,\,
{\rm wt}(v_0)={\rm wt}(\emptyset)=0.
\end{eqnarray*}
\def\bv#1{\fsquare(5mm,#1)}
The actions of $e_i$ and $f_i$ on these basis vectors
are given as follows:
\begin{eqnarray*}
&&\hspace{-30pt}
f_0\left(v_0,v_{\ovl 3}, v_{\ovl 2}, v_{\ovl 1},\emptyset\right)
=\left(v_1,v_2,v_3,\emptyset +\frac{1}{2}v_0,\frac{3}{2}v_1
\right),\\
&&\hspace{-30pt}f_1\left(v_1,v_3,v_0,v_{\ovl 2}, 
\emptyset\right)
=\left(v_2,v_0,2v_{\ovl 3}, v_{\ovl 1}\right),\\
&&\hspace{-30pt}f_2\left(v_2,v_{\ovl 3}\right)=
\left(v_3,v_{\ovl 2}\right),\\
&&\hspace{-30pt}e_0\left(v_1,v_2,v_3,v_0,\emptyset\right)=
\left(\emptyset+\frac{1}{2}v_0,v_{\ovl 3}, v_{\ovl 2},
v_{\ovl 1},\frac{3}{2}v_{\ovl 1}\right),\\
&&\hspace{-30pt}e_1\left(v_2,v_0,v_{\ovl 3},v_{\ovl 1}
\right)=
\left(v_1,2v_3,v_0,v_{\ovl 2}\right),\\
&&\hspace{-30pt}e_2\left(v_3,v_{\ovl 2}\right)=
\left(v_2,v_{\ovl 3}\right),
\end{eqnarray*}
%% kono sitano bun, iru?
where we give non-trivial actions only.

%%%%%%%%%%%%%% Section  %%%%%%%%%%%%%%%%
%\renewcommand{\thesection}{\arabic{section}}
\subsection{Affine Geometric Crystal $\cV_1(\TY(D,3,4))$ in $W(\varpi_1)$ }

Let us review the construction of the 
affine geometric crystal $\cV(\TY(D,3,4))$ in $W(\varpi_1)$ 
following \cite{IN}.

For $\xi\in (\frt^*_{\rm cl})_0$, let $t(\xi)$ be the 
translation as in \cite[Sect 4]{K0} and $\wtil\varpi_i$ 
as in \cite{K1}, indeed, 
$\wtil\varpi_i:=\max(1,\frac{2}{(\al_i,\al_i)})\varpi_i$.
Then we have 
\begin{eqnarray*}
&& t(\wtil\varpi_1)=s_0s_1s_2s_1s_2s_1=:w_1,\\
&& t(\text{wt}(v_{\ovl 2}))=s_2s_1s_2s_1s_0s_1=:w_2,
\end{eqnarray*}
Associated with these Weyl group elements $w_1$ and $w_2$,
we define algebraic varieties $\cV_1=\cV_1(\TY(D,3,4))$ and 
$\cV_2=\cV_2(\TY(D,3,4))\subset W(\varpi_1)$ respectively:
\begin{eqnarray*}
&&\hspace{-30pt}\cV_1:=\{V_1(x)
:=Y_0(x_0)Y_1(x_1)Y_2(x_2)Y_1(x_3)Y_2(x_4)Y_1(x_5)
v_1\,\,\vert\,\,x_i\in\bbC^\times,(0\leq i\leq 5)\},\\
&&\hspace{-30pt}\cV_2:=\{V_2(y):=
Y_2(y_2)Y_1(y_1)Y_2(y_4)Y_1(y_3)Y_0(y_0)Y_1(y_5)
v_{\ovl 2}\,\,\vert\,\,y_i\in\bbC^\times,
(0\leq i\leq 5)\}.
\end{eqnarray*}
Owing to the explicit forms of $f_i$'s on $W(\varpi_1)$
as above, we have $f_0^3=0$, $f_1^3=0$ and $f_2^2=0$ 
and then 
\[
Y_i(c)=(1+\frac{f_i}{c}+\frac{f_i^2}{2c^2})\al_i^\vee(c)
\,\,(i=0,1),\q
Y_2(c)=(1+\frac{f_2}{c})\al_2^\vee(c).
\]
We get explicit forms of $V_1(x)\in\cV_1$ 
and $V_2(y)\in\cV_2$ as in \cite{N3}:
\begin{eqnarray*}
&&V_1(x)=\sum_{1\leq i\leq 3}\left(X_iv_{i}+X_{\ovl i}
v_{\ovl i}\right)+X_{0}v_0+X_\emptyset
\emptyset,\\
&&V_2(y)=\sum_{1\leq i\leq 3}\left(Y_iv_{i}+Y_{\ovl i}
v_{\ovl i}\right)+Y_{0}v_0+Y_\emptyset
\emptyset.
\end{eqnarray*}
where the rational functions $X_i$'s and $Y_i$'s are all
positive in $(x_0,\cd,x_5)$ and 
$(y_0,\cd,y_5)$ respectively (as for their explicit forms, see \cite{IN})
and for any $x$ there exist a unique rational function $a(x)$ and $y$
such that $V_2(y)=a(x)V_1(x)$. Using this result, 
we get the positive birational isomorphism 
$\ovl\sigma:\cV_1\longrightarrow \cV_2$ ($V_1(x)\mapsto V_2(y)$)
and we know that its inverse $\ovl\sigma^{-1}$ is also positive.
The actions of $e_0^c$ on $V_2(y)$ 
(respectively $\gamma_0(V_2(y))$ and 
$\vep_0(V_2(y)))$ are 
induced from the ones on 
$Y_2(y_2)Y_1(y_1)Y_2(y_4)Y_1(y_3)Y_0(y_0)Y_1(y_5)$ 
as an element of the geometric crystal $\cV_2$.
We define the action $e_0^c$ on $V_1(x)$ by
\begin{equation}
e_0^cV_1(x)=\ovl\sigma^{-1}\circ e_0^c\circ
\ovl\sigma(V_1(x))).
\label{e0}
\end{equation}
We also define $\gamma_0(V_1(x))$ 
and $\vep_0(V_1(x))$ by 
\begin{equation}
\gamma_0(V_1(x))=\gamma_0(\ovl\sigma(V_1(x))),\qq
\vep_0(V_1(x)):=\vep_0(\ovl\sigma(V_1(x))).
\label{wt0}
\end{equation}
\begin{thm}[\cite{IN}]
Together with (\ref{e0}), 
(\ref{wt0}) on $\cV_1$, we obtain a
positive affine geometric crystal $\chi:=
(\cV_1,\{e_i\}_{i\in I},
\{\gamma_i\}_{i\in I},\{\vep_i\}_{i\in I})$
$(I=\{0,1,2\})$, whose explicit form is as follows:
first we have $e_i^c$, $\gamma_i$ and $\vep_i$
for $i=1,2$ from the formula (\ref{eici}), 
(\ref{vep-i})
and (\ref{gamma-i}).
\begin{eqnarray*}
&&\hspace{-30pt}
e_1^c(V_1(x))=V_1(x_0,\cC_1x_1,x_2,
\cC_3x_3,x_4,\cC_5x_5),\,
e_2^c(V_1(x))=V_1(x_0,x_1,
\cC_2x_2,x_3,\cC_4x_4,x_5),\\
&&\text{where}\\
&&\cC_1=\frac{\frac{c\,{x_0}}{{x_1}} 
+ \frac{{x_0}\,{{x_2}}}{{{x_1}}^2\,{x_3}} 
+\frac{{x_0}\,{{x_2}}\,{{x_4}}}{{{x_1}}^2\,
{{x_3}}^2\,{x_5}}}{\frac{{x_0}}{{x_1}} 
+ \frac{{x_0}\,{{x_2}}}{{{x_1}}^2\,{x_3}} + 
\frac{{x_0}\,{{x_2}}\,
{{x_4}}}{{{x_1}}^2\,{{x_3}}^2\,{x_5}}},\q
\cC_3=\frac{\frac{c\,{x_0}}{{x_1}} 
+ \frac{c\,{x_0}\,{{x_2}}}{{{x_1}}^2\,{x_3}} + 
\frac{{x_0}\,{{x_2}}\,{{x_4}}}{{{x_1}}^2
\,{{x_3}}^2\,{x_5}}}{
\frac{c\,{x_0}}{{x_1}} 
+ \frac{{x_0}\,{{x_2}}}{{{x_1}}^2\,{x_3}} + 
\frac{{x_0}\,{{x_2}}\,{{x_4}}}
{{{x_1}}^2\,{{x_3}}^2\,{x_5}}},\\
&&\cC_5=\frac{c\,\left( \frac{{x_0}}{{x_1}} + 
      \frac{{x_0}\,{{x_2}}}{{{x_1}}^2\,{x_3}} + 
\frac{{x_0}\,{{x_2}}\,{{x_4}}}{{{x_1}}^2
\,{{x_3}}^2\,{x_5}} \right)}{\frac{c\,{x_0}}{{x_1}} 
+ \frac{c\,{x_0}\,{{x_2}}}{{{x_1}}^2\,{x_3}} + 
\frac{{x_0}\,{{x_2}}\,{{x_4}}}
{{{x_1}}^2\,{{x_3}}^2\,{x_5}}},\,
\cC_2=\frac{\frac{c\,{x_1}^3}{{x_2}} 
+ \frac{{x_1}^3\,{x_3}^3}{{{x_2}}^2\,{x_4}}}
{\frac{{x_1}^3}{{x_2}} + \frac{{x_1}^3\,{x_3}^3}{{{x_2}}^2\,{x_4}}},\,
\cC_4=\frac{c\,\left( \frac{{x_1}^3}{{x_2}} + 
\frac{{x_1}^3\,{x_3}^3}{{{x_2}}^2\,{x_4}} \right) }{\frac{c\,{x_1}^3}
{{x_2}} + \frac{{x_1}^3\,{x_3}^3}{{{x_2}}^2\,{x_4}}},\\
&&\vep_1(V_1(x))={\frac{{x_0}}{{x_1}} 
+ \frac{{x_0}\,{{x_2}}}{{{x_1}}^2\,{x_3}} + 
\frac{{x_0}\,{{x_2}}\,
{{x_4}}}{{{x_1}}^2\,{{x_3}}^2\,{x_5}}},\q
\vep_2(V_1(x))={\frac{{x_1}^3}{{x_2}} 
+ \frac{{x_1}^3\,{x_3}^3}{{{x_2}}^2\,{x_4}}},\\
&&\gamma_1(V_1(x))=\frac{x_1^2x_3^2x_5^2}{x_0x_2x_4},\q
\gamma_2(V_1(x))=\frac{x_2^2x_4^2}{x_1^3x_3^3x_5^3}.
\end{eqnarray*}%%2008.10.21
We also have $e_0^c$, $\vep_0$ and $\gamma_0$ on $V_1(x)$
as:
\begin{eqnarray*}
&&e_0^c(V_1(x))=V_1(\frac{D}{c\cdot E}x_0,\frac{F}{c\cdot E}x_1,
\frac{G}{c^3\cdot E^3}x_2,\frac{D\cdot H}{c^2\cdot E\cdot F}x_3,
\frac{D^3}{c^3\cdot G}x_4,\frac{D}{c\cdot H}x_5),\\
&&\vep_0(V_1(x))=\frac{E}{{{x_0}}^3\,{{x_2}}\,{x_3}},\qq
\gamma_0(V_1(x))=\frac{x_0^2}{x_1x_3x_5},
\end{eqnarray*}
\text{where}\\
\begin{eqnarray*}
&&\hspace{-20pt}
D=c^2\,{{x_0}}^2\,{x_2}\,{x_3} + 
 {x_1}\,{x_2}\,{{x_3}}^2\,{x_5} + 
  c\,{x_0}\,\left( {x_1}\,{{x_3}}^3 + 
 {x_2}\,\left( {{x_3}}^2 + {x_1}\,{x_4} + 
 {x_1}\,{x_3}\,{x_5} \right)  \right)
,\\
&&\hspace{-20pt}
E={{x_0}}^2\,{x_2}\,{x_3} + 
 {x_1}\,{x_2}\,{{x_3}}^2\,{x_5} + 
 {x_0}\,\left( {x_1}\,{{x_3}}^3 + 
 {x_2}\,\left( {{x_3}}^2 + {x_1}\,{x_4} + 
 {x_1}\,{x_3}\,{x_5} \right)  \right),\\
&&\hspace{-20pt}
F= {x_2}\,{{x_3}}^2\,
 \left( {x_0} + {x_1}\,{x_5} \right)  + 
 c\,{x_0}\,\left( {x_0}\,{x_2}\,{x_3} + 
 {x_1}\,\left( {{x_3}}^3 + {x_2}\,{x_4} + 
 {x_2}\,{x_3}\,{x_5} \right)  \right),
\end{eqnarray*}
\begin{eqnarray*}
&&\hspace{-20pt}
G= c^3\,{{x_0}}^6\,{{x_2}}^3\,{{x_3}}^3 + 
      3\,c^2\,{{x_0}}^5\,{{x_2}}^3\,{{x_3}}^4 + 
      3\,c^2\,{{x_0}}^5\,{x_1}\,{{x_2}}^2\,{{x_3}}^5 + 
      3\,c\,{{x_0}}^4\,{{x_2}}^3\,{{x_3}}^5\\
&& + 6\,c\,{{x_0}}^4\,{x_1}\,{{x_2}}^2\,{{x_3}}^6
 +  {{x_0}}^3\,{{x_2}}^3\,{{x_3}}^6 + 
      3\,c\,{{x_0}}^4\,{{x_1}}^2\,{x_2}\,{{x_3}}^7 + 
      3\,{{x_0}}^3\,{x_1}\,{{x_2}}^2\,{{x_3}}^7 \\
&&+  3\,{{x_0}}^3\,{{x_1}}^2\,{x_2}\,{{x_3}}^8
    +  {{x_0}}^3\,{{x_1}}^3\,{{x_3}}^9 + 
3\,c^3\,{{x_0}}^5\,{x_1}\,{{x_2}}^3\,{{x_3}}^2\,{x_4}
+ 6\,c^2\,{{x_0}}^4\,{x_1}\,
  {{x_2}}^3\,{{x_3}}^3\,{x_4}\\
&& + 
3\,c\,{{x_0}}^4\,{{x_1}}^2\,{{x_2}}^2\,{{x_3}}^4\,{x_4}
+3\,c^3\,{{x_0}}^4\,{{x_1}}^2\,
{{x_2}}^2\,{{x_3}}^4\,{x_4} 
+3\,c\,{{x_0}}^3\,{x_1}\,{{x_2}}^3\,{{x_3}}^4\,{x_4}\\
&& + 
3\,{{x_0}}^3\,{{x_1}}^2\,{{x_2}}^2\,{{x_3}}^5\,{x_4}
+3\,c^2\,{{x_0}}^3\,{{x_1}}^2\,
{{x_2}}^2\,{{x_3}}^5\,{x_4}
+ 2\,{{x_0}}^3\,{{x_1}}^3\,{x_2}\,{{x_3}}^6\,{x_4}\\
&& + 
  c^3\,{{x_0}}^3\,{{x_1}}^3\,{x_2}\,{{x_3}}^6\,{x_4}
+ 3\,c^3\,{{x_0}}^4\,{{x_1}}^2\,
{{x_2}}^3\,{x_3}\,{{x_4}}^2 + 3\,c^2\,{{x_0}}^3\,{{x_1}}^2\,{{x_2}}^3\,{{x_3}}^2\,{{x_4}}^2\\
&& + 
 {{x_0}}^3\,{{x_1}}^3\,{{x_2}}^2\,{{x_3}}^3\,{{x_4}}^2 + 
2\,c^3\,{{x_0}}^3\,{{x_1}}^3\,{{x_2}}^2\,
{{x_3}}^3\,{{x_4}}^2
 +  c^3\,{{x_0}}^3\,{{x_1}}^3\,{{x_2}}^3\,{{x_4}}^3\\
&&+3\,c^3\,{{x_0}}^5\,{x_1}\,{{x_2}}^3\,{{x_3}}^3\,{x_5} + 
   9\,c^2\,{{x_0}}^4\,{x_1}\,{{x_2}}^3\,{{x_3}}^4\,{x_5} 
 6\,c^2\,{{x_0}}^4\,{{x_1}}^2\,{{x_2}}^2\,
{{x_3}}^5\,{x_5}\\
&& + 9\,c\,{{x_0}}^3\,{x_1}\,{{x_2}}^3\,{{x_3}}^5\,{x_5} + 
12\,c\,{{x_0}}^3\,{{x_1}}^2\,{{x_2}}^2\,{{x_3}}^6\,{x_5} 
+  3\,{{x_0}}^2\,{x_1}\,{{x_2}}^3\,{{x_3}}^6\,{x_5}\\
&& + 3\,c\,{{x_0}}^3\,{{x_1}}^3\,{x_2}\,{{x_3}}^7\,{x_5} + 
  6\,{{x_0}}^2\,{{x_1}}^2\,{{x_2}}^2\,{{x_3}}^7\,{x_5} + 
      3\,{{x_0}}^2\,{{x_1}}^3\,{x_2}\,{{x_3}}^8\,{x_5} \\
&&+6\,c^3\,{{x_0}}^4\,{{x_1}}^2\,{{x_2}}^3\,
{{x_3}}^2\,{x_4}\, {x_5} 
+ 12\,c^2\,{{x_0}}^3\,{{x_1}}^2\,{{x_2}}^3\,{{x_3}}^3\,
  {x_4}\,{x_5} + 3\,c\,{{x_0}}^3\,{{x_1}}^3\,{{x_2}}^2\,
       {{x_3}}^4\,{x_4}\,{x_5} \\
&&+3\,c^3\,{{x_0}}^3\,{{x_1}}^3\,{{x_2}}^2\,
{{x_3}}^4\,{x_4}\,{x_5} 
+ 6\,c\,{{x_0}}^2\,{{x_1}}^2\,{{x_2}}^3\,{{x_3}}^4\,
{x_4}\,{x_5} + 3\,{{x_0}}^2\,{{x_1}}^3\,{{x_2}}^2\,
       {{x_3}}^5\,{x_4}\,{x_5}
\end{eqnarray*}
\begin{eqnarray*}
&& + 3\,c^2\,{{x_0}}^2\,{{x_1}}^3\,
{{x_2}}^2\,{{x_3}}^5\,{x_4}\,
       {x_5} + 3\,c^3\,{{x_0}}^3\,
{{x_1}}^3\,{{x_2}}^3\,{x_3}\,
{{x_4}}^2\,{x_5} + 3\,c^2\,{{x_0}}^2\,
{{x_1}}^3\,{{x_2}}^3\,
{{x_3}}^2\,{{x_4}}^2\,{x_5} \\
&&+  3\,c^3\,{{x_0}}^4\,{{x_1}}^2\,
{{x_2}}^3\,{{x_3}}^3\,{{x_5}}^2 + 
 9\,c^2\,{{x_0}}^3\,{{x_1}}^2\,{{x_2}}^3\,
{{x_3}}^4\,{{x_5}}^2 + 
 3\,c^2\,{{x_0}}^3\,{{x_1}}^3\,
{{x_2}}^2\,{{x_3}}^5\,{{x_5}}^2 \\
&&+ 9\,c\,{{x_0}}^2\,{{x_1}}^2\,{{x_2}}^3\,
{{x_3}}^5\,{{x_5}}^2 + 
 6\,c\,{{x_0}}^2\,{{x_1}}^3\,{{x_2}}^2\,{{x_3}}^6\,{{x_5}}^2 + 
3\,{x_0}\,{{x_1}}^2\,{{x_2}}^3\,{{x_3}}^6\,{{x_5}}^2 \\
&&+ 3\,{x_0}\,{{x_1}}^3\,{{x_2}}^2\,{{x_3}}^7\,{{x_5}}^2 
+3\,c^3\,{{x_0}}^3\,{{x_1}}^3\,{{x_2}}^3\,
{{x_3}}^2\,{x_4}\,{{x_5}}^2 + 
6\,c^2\,{{x_0}}^2\,{{x_1}}^3\,{{x_2}}^3\,{{x_3}}^3\,
   {x_4}\,{{x_5}}^2\\
&& + 3\,c\,{x_0}\,{{x_1}}^3\,{{x_2}}^3\,
       {{x_3}}^4\,{x_4}\,{{x_5}}^2 
+  c^3\,{{x_0}}^3\,{{x_1}}^3\,{{x_2}}^3\,{{x_3}}^3\,
{{x_5}}^3 +   3\,c^2\,{{x_0}}^2\,{{x_1}}^3\,
{{x_2}}^3\,{{x_3}}^4\,{{x_5}}^3 \\
&&+3\,c\,{x_0}\,{{x_1}}^3\,{{x_2}}^3\,
{{x_3}}^5\,{{x_5}}^3 + 
 {{x_1}}^3\,{{x_2}}^3\,{{x_3}}^6\,{{x_5}}^3,\\
&&\hspace{-20pt}
H= c\,{{x_0}}^2\,{x_2}\,{x_3} + 
 {x_0}\,{x_2}\,{{x_3}}^2 + {x_0}\,{x_1}\,{{x_3}}^3 + 
 {x_0}\,{x_1}\,{x_2}\,{x_4} + 
 c\,{x_0}\,{x_1}\,{x_2}\,{x_3}\,{x_5} + 
 {x_1}\,{x_2}\,{{x_3}}^2\,{x_5}.
\end{eqnarray*}
\end{thm}

%%%%%%%%%%%%%% Section  %%%%%%%%%%%%%%%%
\renewcommand{\thesection}{\arabic{section}}
\section{Ultra-discretization}
\setcounter{equation}{0}
\renewcommand{\theequation}{\thesection.\arabic{equation}}

We denote the positive structure on $\chi$ as in the 
previous section by 
$\theta:T'\seteq(\bbC^\times)^6 \longrightarrow \cV_1$
($x\mapsto V_1(x)$).
Then by Corollary \ref{cor-posi}
we obtain the ultra-discretization 
${\mathcal UD}(\chi,T',\theta)$, 
which is a Kashiwara's crystal. 
Now we show that the conjecture in \cite{IN} 
is correct and it turns 
out to be the following theorem.
\begin{thm}
\label{ultra-d}
The crystal ${\mathcal UD}(\chi,T',\theta)$ as above 
is isomorphic to the crystal $B_\ify$ of 
type $\TY(G,1,2)$ as in Sect.\ref{perf}.
\end{thm}
In order to show the theorem, we shall see
the explicit crystal structure on 
$\cX:={\mathcal UD}(\chi,T',\theta)$.
Note that ${\mathcal UD}(\chi)=\ZZ^6$ as a set .
Here as  for variables in $\cX$, 
we use the same notations $c,x_0,x_1,\cd,x_5$
as for $\chi$.

For $x=(x_0,x_1,\cd,x_5)\in\cX$, it follows from the 
results in the previous section that the functions
$\wt_i$ and $\vep_i$ ($i=0,1,2$) are given as:
\begin{eqnarray*}
&&\wt_0(x)=2x_0-x_1-x_3-x_5,\,\,
\wt_1(x)=2(x_1+x_3+x_5)-x_0-x_2-x_4,\\
&&\wt_2(x)=2(x_2+x_4)-3(x_1-x_3-x_5).
\end{eqnarray*}
Set
\begin{equation}
\begin{array}{l}
\al\seteq 2x_0+x_2+x_3,\q
\beta\seteq x_1+x_2+2x_3+x_5,\q
\gamma\seteq x_0+x_1+3x_3,\,\\
\del\seteq x_0+x_2+2x_3,\q
\epsilon\seteq x_0+x_1+x_2+x_4,\q
\\ 
\phi\seteq x_0+x_1+x_2+x_3+x_5.
\end{array}
\label{alphabeta}
\end{equation}
Then we have
\begin{eqnarray}
&&\vep_0(x)=\max(\al,\beta,\gamma,\del,\epsilon,
\phi)-(3x_0+x_2+x_3),\nn \\
&&\vep_1(x)=\max(x_0-x_1,x_0+x_2-2x_1-x_3,
x_0+x_2+x_4-2x_1-2x_3-x_5),\\
&&\vep_2(x)=\max(3x_1-x_2,3x_1+3x_3-2x_2-x_4).\nn
\end{eqnarray}
Indeed, from the explicit form of $G$ as in the previous
section we have
\begin{eqnarray*}
 &&{\mathcal UD}(G)|_{c=-1}
=\max(-3+3\al,-2+2\al+\del,-2+2\al+\gamma,-1+\al+2\del,
 -1+\al+\gamma+\del,\\
 && 3\del,-1+\al+2\gamma,\gamma+2\del,2\gamma+\del,
 3\gamma,-3+2\al+\epsilon,-2+\al+\del+\epsilon,-1+\al+\gamma+\epsilon,\\
&& -1+2\del+\epsilon,\gamma+\del+\epsilon,2\gamma+\epsilon,-3+\al+2\epsilon,
 -2+\del+2\epsilon,\gamma+2\epsilon,-3+3\epsilon,-3+2\al+\phi,\\
 && -2+\al+\del+\phi,-2+\al+\gamma+\phi,-1+2\del+\phi,-1+\gamma+\del+\phi,\beta+2\del,-1+2\gamma+\phi,\\
 &&\beta+\gamma+\del,\beta+2\gamma,-3+\al+\epsilon+\phi,-2+\del+\epsilon+\phi,
 -1+\gamma+\epsilon+\phi,-1+\beta+\del+\epsilon,\\
 &&\beta+\gamma+\epsilon,-3+2\epsilon+\phi,-2+\beta+2\epsilon,-3+\al+2\phi,
 -2+\del+2\phi,-2+\gamma+2\phi,\\
 &&-1+\al+2\beta,-1+\beta+\gamma+\phi,2\beta+\del,
 2\beta+\gamma,-3+\epsilon+2\phi,-2+\beta+\epsilon+\phi,\\
 &&-1+2\beta+\epsilon,-3+3\phi,-2+\beta+2\phi,-1+2\beta+\phi,3\beta).
\end{eqnarray*}
We simplify this by using the following lemma:
\begin{lem}
\label{convex}
For $m_1,\cd,m_k\in\RR$ and $t_1,\cd,t_k\in
\RR_{\geq0}$ such that $t_1+\cd t_k=1$,
we have 
\[
 \max\left(m_1,\cd,m_k,\sum_{i=1}^kt_im_i\right)
= \max(m_1,\cd,m_k)
\]
\end{lem}
Since we have 
\begin{eqnarray*}\label{3tobun}
&&-2+2\al+\del=\frac{2(-3+3\al)+3\del}{3},\q
-2+2\al+\gamma=\frac{2(-3+3\al)+3\gamma}{3},\\
&&-1+\al+2\del=\frac{2\cdot3\del+(-3+3\al)}{3},\q
-1+\al+\gamma+\del=\frac{(-3+3\al)+3\gamma+3\del}{3},\\
&&-1+\al+2\gamma=\frac{(-3+3\al)+2\cdot3\gamma}{3},\q
\gamma+2\del=\frac{2\cdot 3\del+3\gamma}{3},\q
\q{\rm etc,}
\end{eqnarray*}
by this lemma we get
\begin{eqnarray*}
&&{\mathcal UD}(G)|_{c=-1}
=\max(-3+3\al,3\beta,3\gamma,3\del,-3+3\epsilon,-3+3\phi,
-1+\al+\gamma+\epsilon,\gamma+\del+\epsilon,\\
&&\gamma+2\epsilon,2\gamma+\epsilon,
-1+\gamma+\epsilon+\phi,\beta+\gamma+\epsilon).
\end{eqnarray*}
Next, we describe the actions of $\fit$ $(i=0,1,2)$.
Set $\Xi_j\seteq{\mathcal UD}(\cC_j)|_{c=-1}$ 
($j=1,\cd,5$). Then we have
%%c=-1 wo dainyu? c nomama?
\begin{eqnarray*}
\Xi_1&=&\max(-1+x_0-x_1,x_0+x_2-2x_1-x_3,
x_0+x_2+x_4-2x_1-2x_3-x_5)\\
&&-\max(x_0-x_1,x_0+x_2-2x_1-x_3,
x_0+x_2+x_4-2x_1-2x_3-x_5),\\
\Xi_3&=&\max(-1+x_0-x_1,-1+x_0+x_2-2x_1-x_3,
x_0+x_2+x_4-2x_1-2x_3-x_5)\\
&&-\max(-1+x_0-x_1,x_0+x_2-2x_1-x_3,
x_0+x_2+x_4-2x_1-2x_3-x_5),\\
\Xi_5&=&\max(-1+x_0-x_1,-1+x_0+x_2-2x_1-x_3,
-1+x_0+x_2+x_4-2x_1-2x_3-x_5)\\
&&-\max(-1+x_0-x_1,-1+x_0+x_2-2x_1-x_3,
x_0+x_2+x_4-2x_1-2x_3-x_5),\\
\Xi_2&=&\max(-1+3x_1-x_2,3x_1+3x_3-2x_2-x_4)
-\max(3x_1-x_2,3x_1+3x_3-2x_2-x_4),\\
\Xi_4&=&\max(-1+3x_1-x_2,-1+3x_1+3x_3-2x_2-x_4)\\
&&-\max(-1+3x_1-x_2,3x_1+3x_3-2x_2-x_4).
\end{eqnarray*}
Therefore, for $x\in\cX$ we have
\begin{eqnarray*}
&&\til f_1(x)=(x_0,x_1+\Xi_1,x_2,
x_3+\Xi_3,x_4,x_5+\Xi_5),\\
&&\til f_2(x)=(x_0,x_1,x_2+\Xi_2,x_3,x_4+\Xi_4,x_5).
\end{eqnarray*}
We obtain the action $\eit$ ($i=1,2$)
by setting $c=1$ in ${\mathcal UD}(\cC_i)$.
Finally, we describe the action of $\til f_0$. 
Set
\begin{eqnarray*}
\Psi_0&\seteq&
\max(-2+\al,\beta,-1+\gamma,-1+\del,-1+\epsilon,
-1+\phi)\\
&&-\max(\al,\beta,\gamma,\del,\epsilon,\phi)+1,
\\
\Psi_1&\seteq&
\max(-1+\al,\beta,-1+\gamma,\del,-1+\epsilon,
-1+\phi)\\
&&-\max(\al,\beta,\gamma,\del,\epsilon,\phi)+1,
\\
\Psi_2&\seteq&
\max(-3+3\al,3\beta,3\gamma,3\del,-3+3\epsilon,
-3+3\phi,-1+\al+\gamma+\epsilon,\gamma+\del+\epsilon,\\
&&\gamma+2\epsilon,2\gamma+\epsilon,-1+\gamma+\epsilon+\phi,
\beta+\gamma+\epsilon)\\
&&-3\max(\al,\beta,\gamma,\del,\epsilon,\phi)+3,
\\
\Psi_3&\seteq&
\max(-2+\al,\beta,-1+\gamma,-1+\del,-1+\epsilon,
-1+\phi)\\
&&+\max(-1+\al,\beta,\gamma,\del,\epsilon,-1+\phi)
-\max(\al,\beta,\gamma,\del,\epsilon,\phi)\\
&&-\max(-1+\al,\beta,-1+\gamma,\del,-1+\epsilon,
-1+\phi)+2,\\
\Psi_4&\seteq&
3\max(-2+\al,\beta,-1+\gamma,-1+\del,-1+\epsilon,
-1+\phi)\\
&&-\max(-3+3\al,3\beta,3\gamma,3\del,-3+3\epsilon,-3+3\phi,
-1+\al+\gamma+\epsilon,\gamma+\del+\epsilon,\\
&&\gamma+2\epsilon,
2\gamma+\epsilon,-1+\gamma+\epsilon+\phi,\beta+\gamma+\epsilon)+3,
\\
\Psi_5&\seteq&
\max(-2+\al,\beta,-1+\gamma,-1+\del,-1+\epsilon,
-1+\phi)\\
&&-\max(1+\al,\beta,\gamma,\del,\epsilon,-1+\phi)+1,
\end{eqnarray*}
where $\al,\beta,\cd,\phi$ are as in (\ref{alphabeta}).
Therefore, by the explicit form of $e_0^c$ 
as in the previous section, we have
\begin{equation}
\til f_0(x)=(x_0+\Psi_0,x_1+\Psi_1,x_2+\Psi_2,
x_3+\Psi_3,x_4+\Psi_4,x_5+\Psi_5).
\end{equation}
We have the explicit form of $\til e_0$ by setting $c=1$ in 
${\mathcal UD}({\mathcal C}_i)$.
Now, let us show the theorem.\\
({\sl Proof of Theorem \ref{ultra-d}.})
Define the map 
\[
\begin{array}{cccc}
\Omega\cl&\cX&\longrightarrow& B_\ify,\\
&(x_0,\cd, x_5)&\mapsto& (b_1,b_2,b_3,
\ovl b_3,\ovl b_2,\ovl b_1),
\end{array}
\]
by 
\[
 b_1=x_5,\,\, b_2=\frac{1}{3}x_4-x_5,\,\,
b_3=x_3-\frac{2}{3}x_4,\,\, \ovl b_3=\frac{2}{3}x_2-x_3,\,\,
\ovl b_2=x_1-\frac{1}{3}x_2,\,\,
\ovl b_1=x_0-x_1,
\]
and $\Omega^{-1}$ is given by 
\begin{eqnarray*}
&&x_0=b_1+b_2+\frac{b_3+\ovl b_3}{2}+\ovl b_2+\ovl b_1,
\q
x_1=b_1+b_2+\frac{b_3+\ovl b_3}{2}+\ovl b_2,\\
&&
x_2=3b_1+3b_2+\frac{3(b_3+\ovl b_3)}{2},\,\,
x_3=2b_1+2b_2+b_3,\,\,
x_4=3b_1+3b_2,\,\,
x_5=b_1,
\end{eqnarray*}
which means that $\Omega$ is bijective.
Here note that 
$\frac{3(b_3+\ovl b_3)}{2}\in \bbZ$ by the 
definition of $B_\ify$ as in Sect.4.
We shall show that 
$\Omega$ is commutative with actions of $\fit$
and preserves the functions $\wt_i$ and $\vep_i$,
that is, 
\[
 \fit(\Omega(x))=\Omega(\fit x),\q
\wt_i(\Omega(x))=\wt_i(x),\q
\vep_i(\Omega(x))=\vep_i(x)
\q(i=0,1,2),
\]
Indeed, the commutativity 
$\eit(\Omega(x))=\Omega(\eit x)$ is shown by a similar way.
First, let us check $\wt_i$: Set $b=\Omega(x)$ and 
let $(z_1,z_2,z_3,z_4)$ be as in (\ref{z1-4}).
By the explicit forms of $\wt_i$ on 
$\cX$ and $B_\ify$,
we have 
\begin{eqnarray*}
&&\hspace{-10pt}
\wt_0(\Omega(x))=\vp_0(\Omega(x))
-\vep_0(\Omega(x))=2z_1+z_2+z_3+3z_4\q\qq\\
&&=
2(\ovl b_1-b_1)+(\ovl b_2-\ovl b_3)
+(b_3-b_2)+\frac{3}{2}(\ovl b_3-b_3)
=2(\ovl b_1-b_1)+\ovl b_2-b_2
+\frac{\ovl b_3-b_3}{2}\\
&&=2x_0-x_1-x_3-x_5=\wt_0(x),\\
&&\hspace{-10pt}\wt_1(\Omega(x))
=\vp_1(\Omega(x))-\vep_1(\Omega(x))\\
&&=
b_1+(b_3-b_2+({\bar b}_2-{\bar b}_3)_+)_+
-({\bar b}_1+({\bar b}_3
-{\bar b}_2+(b_2-b_3)_+)_+)\\
&&=b_1-\ovl b_1-b_2+\ovl b_2+b_3-\ovl b_3
=2(x_1+x_3+x_5)-x_0-x_2-x_4=\wt_1(x),\\
&&\hspace{-10pt}
\wt_2(\Omega(x))
=\vp_2(\Omega(x))-\vep_2(\Omega(x))
=
3b_2+\frac{3}{2}({\bar b}_3-b_3)_+-
3{\bar b}_2-\frac{3}{2}(b_3-{\bar b}_3)_+\\
&&=3b_2-3\ovl b_2+\frac{3}{2}(\ovl b_3-b_3)
=2(x_2+x_4)-3(x_1+x_3+x_5)=\wt_2(x).
\end{eqnarray*}
Next, we shall check $\vep_i$:
\begin{eqnarray*}
&&\vep_1(\Omega(x))
={\bar b}_1
+({\bar b}_3-{\bar b}_2+(b_2-b_3)_+)_+\\
&&\,\,=\max({\bar b}_1,{\bar b}_1+
{\bar b}_3-{\bar b}_2,
{\bar b}_1+
{\bar b}_3-{\bar b}_2+b_2-b_3)\\
&&\,\,
=\max(x_0-x_1,x_0-2x_1+x_2-x_3,
x_0-2x_1+x_2-2x_3+x_4-x_5)=\vep_1(x),\\
&&\vep_2(\Omega(x))
=3{\bar b}_2+\frac{3}{2}(b_3-{\bar b}_3)_+
=\max(3{\bar b}_2,
3{\bar b}_2+\frac{3}{2}(b_3-{\bar b}_3))\\
&&\,\,
=\max(3x_1-x_2,3x_1-2x_2+3x_3-x_4)
=\vep_2(x).
\end{eqnarray*}
Here let us see $\vep_0$:
\begin{eqnarray*}
&&\vep_0(\Omega(x))=
-s(b)+\max A-(2z_1+z_2+z_3+3z_4)\\
&&=-x_0+\max(0,z_1,z_1+z_2,z_1+z_2+3z_4,
z_1+z_2+z_3+3z_4,2z_1+z_2+z_3+3z_4)-(\al-\beta)
\\&&
=-x_0+\max
(-2x_0+x_1+x_3+x_5,-x_0+x_3,-x_0+x_1-x_2+2x_3,
\\&&
\qq\qq \qq\qq -x_0+x_1-x_3+x_4,-x_0+x_1+x_5,0)\\
&&=
-(3x_0+x_2+x_3)+\max(x_1+x_2+2x_3+x_5,
x_0+x_2+2x_3,x_0+x_1+3x_3,\\
&&\qq x_0+x_1+x_2+x_4,
x_0+x_1+x_2+x_3+x_5,2x_0+x_2+x_3)\\
&&=-(3x_0+x_2+x_3)
+\max(\beta,\del,\gamma,\epsilon,\phi,\al).
\end{eqnarray*}
On the other hand, we have
\[
\vep_0(x)=-(3x_0+x_2+x_3)
+\max(\al,\beta,\gamma,\del,\epsilon,
\phi).
\]
which shows
$\vep_0(\Omega(x))=\vep_0(x)$.

Let us show
 $\fit(\Omega(x))=\Omega(\fit(x))$
($x\in\cX,\,i=0,1,2$).
As for $\til f_1$, set 
\[
 A=x_0-x_1,\,\, B=x_0+x_2-2x_1-x_3,\,\,
C=x_0+x_2+x_4-2x_1-2x_3-x_5.
\]
Then we obtain
$\Xi_1=\max(A-1,B,C)-\max(A,B,C),\,\,
\Xi_3=\max(A-1,B-1,C)-\max(A-1,B,C),\,\,
\Xi_5=\max(A-1,B-1,C-1)-\max(A-1,B-1,C).$
%% > < nomuki tyuui!! \geq toka = nohairuhairanai mo
Therefore, we have
\begin{eqnarray*}
&&\Xi_1=-1,\,\,\Xi_3=0,\,\,\Xi_5=0,\,\,
\text{if}\,\,A> B,C\\
&&\Xi_1=0,\,\,\Xi_3=-1,\,\,\Xi_5=0,\,\,
\text{if}\,\,A\leq B> C\\
&&\Xi_1=0,\,\,\Xi_3=0,\,\,\Xi_5=-1,\,\,
\text{if}\,\,A, B\leq C, 
\end{eqnarray*}
which implies
\begin{eqnarray*}
\til f_1(x)=\begin{cases}
(x_0,x_1-1,x_2,\cd,x_5)&\text{if }
A> B,C\\
(x_0,\cd,x_3-1,x_4,x_5)&\text{if }
A\leq B> C\\
(x_0,\cd,x_4,x_5-1)&\text{if }
A, B\leq C
\end{cases}
\end{eqnarray*}
Since $A=\ovl b_1$, 
$B=\ovl b_1+\ovl b_3-\ovl b_2$ and 
$C=\ovl b_1+\ovl b_3-\ovl b_2+b_2-b_3$,
we get ($b=\Omega(x)$)
\begin{eqnarray*}
\Omega({\tilde f}_1(x))=&
\begin{cases}
(\ldots,{\bar b}_2 -1,{\bar b}_1 +1) 
& \text{if ${\bar b}_2 -{\bar b}_3 > (b_2 -b_3)_+$}, 
\\  
(\ldots,b_3 -1,{\bar b}_3 +1,\ldots) 
& \text{if ${\bar b}_2 -{\bar b}_3\leq 0 < b_3 -b_2$}, 
\\ 
(b_1 -1,b_2 +1,\ldots) 
& \text{if $({\bar b}_2 -{\bar b}_3)_+ \leq b_2 -b_3$},
\end{cases}
\end{eqnarray*}
which is the same as the action of 
$\til f_1$ on $b=\Omega(x)$ as in Sect.4.
Hence, we have 
$\Omega(\til f_1(x))=\til f_1(\Omega(x))$.

Let us see 
$\Omega(\til f_2(x))=\til f_2(\Omega(x))$.
Set
\[
 L=3x_1-x_2, \q M:=3x_1+3x_3-2x_2-x_4.
\]
Then $\Xi_2=\max(-1+L,M)-\max(L,M)$ and 
$\Xi_4=\max(-1+L,-1+M)-\max(-1+L,M)$. Thus,
one has
\begin{eqnarray*}
&&\Xi_2=-1,\q\Xi_4=0 \q\text{if }L>  M,\\
&&\Xi_2=0,\q\Xi_4=-1 \q\text{if }L\leq M,
\end{eqnarray*}
which means
\[
 \til f_2(x)=\begin{cases}
(x_0,x_1,x_2-1,x_3,x_4,x_5)&\text{if }L> M,\\
(x_0,x_1,x_2,x_3,x_4-1,x_5)&\text{if }L\leq M.
\end{cases}
\]
Since $L-M=x_2-3x_3+x_4=\frac{3(\ovl b_3-b_3)}{2}$,
one gets
\[
\Omega( \til f_2(x))=
\begin{cases}
(\ldots,{\bar b}_3 -\frac{2}{3},{\bar b}_2 +\frac{1}{3},\ldots) 
& \text{if ${\bar b}_3 > b_3$}, 
\\
(\ldots,b_2 -\frac{1}{3},b_3 +\frac{2}{3},\ldots) 
& \text{if ${\bar b}_3 \leq b_3$},
\end{cases}
\]
where $b=\Omega(x)$. 
This action coincides with the one of
 $\til f_2$ on $b\in B_\ify$ as in Sect4.
Therefore, we get 
$\Omega(\til f_2(x))=\til f_2(\Omega(x))$.

Finally, we shall check 
$\til f_0(\Omega(x))=\Omega(\til f_0(x))$.
For the purpose, we shall estimate the values
$\Psi_0,\cd,\Psi_5$ explicitly.

First, the following cases are investigated:
\begin{eqnarray*}
&({\rm f}1)&\beta\geq \gamma,\del,\epsilon,
\phi,\,\,\, \phi\geq\al,\,\, \del\geq\al\\
&({\rm f}2)&\beta< \del\geq \al,\gamma,\epsilon,
\,\,\,\al>\phi,\,\,\beta\geq\phi\\
&({\rm f}3)&\beta,\del< \gamma\geq \al,\epsilon,
\phi\\
&({\rm f}4)&\beta,\del < \epsilon \geq\al,\phi,\,\,\epsilon=\gamma+1\\
&({\rm f}4')&\beta,\del < \epsilon \geq \al,\phi,\,\,\epsilon=\gamma+2\\
&({\rm f}4'')&\beta,\del < \epsilon \geq \al,\phi,\,\,\epsilon > \gamma+2\\
&({\rm f}5)&\beta,\gamma,\epsilon < \phi\geq \al,\,\,\,\al>\del,\,\,\beta\geq\del\\
&({\rm f}6)&\al > \gamma,\del,\epsilon,\phi,\,\,\,\del,\phi>\beta.
\end{eqnarray*}
It is  easy to see that each of these
 conditions are equivalent to the conditions 
$(F_1)$-$(F_6)$ in Sect.4, more precisely, we have
$({\rm f}i)\Leftrightarrow\,\, (F_i)$ ($i=1,2,3,5,6$), 
$({\rm f}4)\Leftrightarrow\ (F_4)$ and $z_4=\frac{1}{3}$,
$({\rm f}4')\Leftrightarrow\ (F_4)$ and $z_4=\frac{2}{3}$
and $({\rm f}4'')\Leftrightarrow\ (F_4)$ and $z_4\ne\frac{1}{3}, \frac{2}{3}$,
and that (f1)--(f6) cover all cases and they 
have no intersection.

%% kokokara naositenai
Let us show (f1)$\Leftrightarrow\,(F_1)$:
the condition (f1) means
$\beta-\gamma=-(z_1+z_2)\geq 0$, 
$\beta-\del=-z_1\geq 0$, 
$\beta-\epsilon=-(z_1+z_2+3z_4)\geq 0$ and
$\beta-\phi=-(z_1+z_2+z_3+3z_4)\geq 0$,
which is equivalent to the 
condition $z_1+z_2\leq0$, $z_1\leq0$, 
$z_1+z_2+3z_4\leq0$ and 
$z_1+z_2+z_3+3z_4\leq0$. (Note that $\phi-\al=\beta-\del,
\del-\al=\beta-\phi$) This is just 
the condition $(F_1)$.
Other cases $i=2,3,5,6$ are shown similarly.
Next, let us see the cases $({\rm f}4)$, $({\rm f}4')$ and 
$({\rm f}4'')$. Indeed, 
\[
\epsilon-\gamma=x_2-3x_3+x_4=\frac{3}{2}(\ovl b_3-b_3)=3z_4.
\]
Thus, we can easily get that 
$({\rm f}4)\Leftrightarrow\ (F_4)$ and $z_4=\frac{1}{3}$,
$({\rm f}4')\Leftrightarrow\ (F_4)$ and $z_4=\frac{2}{3}$.
and $({\rm f}4'')\Leftrightarrow\ (F_4)$ and $z_4\ne\frac{1}{3}, \frac{2}{3}$.

Under the condition (f1)
($\Leftrightarrow\,(F_1)$),
we have 
\[
 \Psi_0=\Psi_1=\Psi_5=1,\Psi_2=\Psi_4=3,\q
\Psi_3=2,
\]
which means
$\til f_0(x)=(x_0+1,x_1+1,x_2+3,x_3+2,x_4+3,
x_5+1)$. Thus, we have
\[
 \Omega(\til f_0(x))=(b_1+1,b_2,\cd,\ovl b_1),
\]
which coincides with the action of 
$\til f_0$ under $(F_1)$ in Sect.4.
Similarly, we have
\[
\begin{array}{ccc}
{\rm(f2)}&\Rightarrow&
(\Psi_0,\Psi_1,\Psi_2,\Psi_3,\Psi_4,\Psi_5)
=(0,1,3,1,0,0)\\
&\Rightarrow
&\til f_0(x)=(x_0,x_1+1,x_2+3,x_3+1,x_4,x_5),\\
&\Rightarrow&
\Omega(\til f_0(x))
=(b_1,b_2,b_3+1,\ovl b_3+1,\ovl b_2,\ovl b_1-1),
\end{array}
\]
which coincides with the action of $\til f_0$
under $(F_2)$ in Sect.4.
\[
\begin{array}{ccc}
{\rm(f3)}&\Rightarrow&
(\Psi_0,\Psi_1,\Psi_2,\Psi_3,\Psi_4,\Psi_5)
=(0,0,3,2,0,0)\\
&\Rightarrow
&\til f_0(x)=(x_0,x_1,x_2+3,x_3+2,x_4,x_5),\\
&\Rightarrow&
\Omega(\til f_0(x))
=(b_1,b_2,b_3+2,\ovl b_3,\ovl b_2-1,\ovl b_1),
\end{array}
\]
which coincides with the action of $\til f_0$
under $(F_3)$ in Sect.4.
%% kokomade konosita mada
\[
\begin{array}{ccc}
\text{(f4)}&\Rightarrow&
(\Psi_0,\Psi_1,\Psi_2,\Psi_3,\Psi_4,\Psi_5)
=(0,0,2,2,1,0)\\
&\Rightarrow
&\til f_0(x)=(x_0,x_1,x_2+2,x_3+2,x_4+1,x_5),\\
&\Rightarrow&
\Omega(\til f_0(x))
=(b_1,b_2+\frac{1}{3},b_3+\frac{4}{3},\ovl b_3-\frac{2}{3},\ovl b_2-\frac{2}{3},\ovl b_1),
\end{array}
\]
which coincides with the action of $\til f_0$
under $(F_4)$ and $z_4=\frac{1}{3}$ in Sect.4.
\[
\begin{array}{ccc}
\text{(f4$'$)}&\Rightarrow&
(\Psi_0,\Psi_1,\Psi_2,\Psi_3,\Psi_4,\Psi_5)
=(0,0,1,2,2,0)\\
&\Rightarrow
&\til f_0(x)=(x_0,x_1,x_2+1,x_3+2,x_4+2,x_5),\\
&\Rightarrow&
\Omega(\til f_0(x))
=(b_1,b_2+\frac{2}{3},b_3+\frac{2}{3},\ovl b_3-\frac{4}{3},\ovl b_2-\frac{1}{3},\ovl b_1),
\end{array}
\]
which coincides with the action of $\til f_0$
under $(F_4)$ and $z_4=\frac{2}{3}$ in Sect.4.
\[
\begin{array}{ccc}
\text{(f4$''$)}&\Rightarrow&
(\Psi_0,\Psi_1,\Psi_2,\Psi_3,\Psi_4,\Psi_5)
=(0,0,0,2,3,0)\\
&\Rightarrow
&\til f_0(x)=(x_0,x_1,x_2,x_3+2,x_4+3,x_5),\\
&\Rightarrow&
\Omega(\til f_0(x))
=(b_1,b_2+1,b_3,\ovl b_3-2,\ovl b_2,\ovl b_1),
\end{array}
\]
which coincides with the action of $\til f_0$
under $(F_4)$ and $z_4\ne\frac{1}{3},\frac{2}{3}$ in Sect.4.
%%konoue naosu
\[
\begin{array}{ccc}
{\rm(f5)}&\Rightarrow&
(\Psi_0,\Psi_1,\Psi_2,\Psi_3,\Psi_4,\Psi_5)
=(0,0,0,1,3,1)\\
&\Rightarrow
&\til f_0(x)=(x_0,x_1,x_2,x_3+1,x_4+3,x_5+1),\\
&\Rightarrow&
\Omega(\til f_0(x))
=(b_1+1,b_2,b_3-1,\ovl b_3-1,\ovl b_2,\ovl b_1),
\end{array}
\]
which coincides with the action of $\til f_0$
under $(F_5)$ in Sect.4.
\[
\begin{array}{ccc}
{\rm(f6)}&\Rightarrow&
(\Psi_0,\Psi_1,\Psi_2,\Psi_3,\Psi_4,\Psi_5)
=(-1,0,0,0,0,0)\\
&\Rightarrow
&\til f_0(x)=(x_0-1,x_1,x_2,x_3,x_4,x_5),\\
&\Rightarrow&
\Omega(\til f_0(x))
=(b_1,b_2,b_3,\ovl b_3,\ovl b_2,\ovl b_1-1),
\end{array}
\]
which coincides with the action of $\til f_0$
under $(F_6)$ in Sect.4.
Now, we have 
$\Omega(\til f_0(x))=\til f_0(\Omega(x))$.
Therefore, 
the proof of Theorem \ref{ultra-d} has been
completed.\qed

\bibliographystyle{amsalpha}

\end{document}